    \documentclass[12pt]{article}
    \usepackage{amsmath,amssymb,amsthm,amsfonts}
    \usepackage{paralist}
    \usepackage{subcaption}
\usepackage{graphicx}
    \newtheorem{definition}{Definition}[section]
    \usepackage{graphics} 
    \usepackage{epsfig}
    \newcommand{\sech}{\operatorname{sech}}
    \usepackage{graphicx}
    \usepackage{caption} 
    \usepackage{epstopdf}
    \usepackage[colorlinks=true]{hyperref}
    \usepackage{float}
    \usepackage{pdfpages}
    \usepackage{amsfonts}
    \usepackage{amsmath,amsthm}
    \usepackage{multirow}

    \usepackage{amsmath}
    \usepackage{txfonts}
    \usepackage{enumerate}
    \usepackage{amsmath,bm}
    \usepackage[numbers,sort&compress]{natbib}
    \usepackage[font=scriptsize,labelfont=bf,width=0.95\textwidth]{caption}
    
    \numberwithin{equation}{section}
    \newcommand{\R}{{\mathbb R}}

    \newcommand{\C}{{\mathbb C}}
    
    \newcommand{\be}{\begin{equation}}
    \newcommand{\ee}{\end{equation}}

    \newcommand{\ben}{\begin{eqnarray*}}
    \newcommand{\enn}{\end{eqnarray*}}
    
    \newtheorem{theorem}{\textbf Theorem}[section]
    \newtheorem{lemma}{\textbf Lemma}[section]

    \newtheorem{corollary}{Corollary}

    \numberwithin{equation}{section}
    \newtheorem{remark}{Remark}[section]
    	\usepackage{microtype}

\begin{document}
    
    \title{\textbf{A Complete Classification of the One-Dimensional Nonlocal Eigenvalue Problem}}
    \author{Ziyi Chi
    	\thanks{Department of Mathematics,
    		The Chinese University of Hong Kong,
    		Shatin, New Territories, Hong Kong;
    		email: 1155250023@link.cuhk.edu.hk},~
    	Fanze Kong \thanks{Department of Applied Mathematics, University of Washington, Seattle, 98195, WA, USA, email address: fzkong@uw.edu}~  
    	and Juncheng Wei   \thanks{Department of Mathematics, The Chinese University of Hong Kong, Shatin, New Territories, Hong Kong; \texttt{wei@math.cuhk.edu.hk}.}}
    \date{\today}
    \maketitle

    \vspace{-0.25in}
    
    \begin{abstract}
    We give a complete classification of the spectral stability of the following one-dimensional nonlocal eigenvalue problem (NLEP)
    $$	\partial_x^2 \phi-\phi+pw^{p-1}\phi -\gamma (p-1) \frac{\int_\R w^{r-1}\phi\,dx}{\int_\R w^r\,dx} \, w^p
    	=\lambda\phi,
    $$
    where $ w$ is the unique positive even solution of $ w''-w+w^p=0$ decays at infinity, $p,r>1$  and $\gamma\in\mathbb C$
with $\Re\gamma\ge0$.  This type of NLEP appears frequently in pattern formations in reaction-diffusion systems.   We show that on the even subspace, strict spectral stability holds {exactly for \(\Re\gamma>1\) when \(r=(p+3)/2\).
   In the \(\gamma\)-plane, when \(r>(p+3)/2\), two conjugate
boundary arcs open to the left from \(\gamma=1\), and the
strictly stable region lies to their right. When \(r<(p+3)/2\),
the arcs open to the right, and the strictly stable region
lies between them. Their common endpoint is \(\gamma=1\)
when \(r\ge(p-1)/2\) and \(\gamma=\gamma_H>1\) when
\(1<r<(p-1)/2\). In each case, the two arcs approach a
common vertical asymptote.} As an application, for real \(\gamma\ge0\), the condition for strict
spectral stability is \(\gamma>1\) when \(r\ge(p-1)/2\) and
\(\gamma>\gamma_H\) when \(1<r<(p-1)/2\), where \(\gamma_H>1\)
is the value at which a conjugate pair of eigenvalues crosses
the imaginary axis.
 The key ingredient in determining how the number of eigenvalues
in the open right half-plane changes at imaginary-axis crossings
is a cumulative estimate for the mixed-sign spectral weights
that arise for certain \(p\) and \(r\) with \(1<p<3\) and \(1<r<2\).

    \end{abstract}

    
    
    \maketitle

    \section{Introduction}
    
    Nonlocal eigenvalue problems (NLEP) arise in the stability analysis
    of localized patterns in reaction--diffusion and chemotaxis
    systems. Applications include the Gierer--Meinhardt,
    Gray--Scott, and Sel'kov models, hypercycle systems and
    models of angiogenesis \cite{WZ01,WW03,WWbook},
    the Schnakenberg and Brusselator models \cite{IWW04,TNW13},
    urban crime models \cite{KWW14}, and Keller--Segel
    systems \cite{SWW05,KKW07,KWW24}.   The mathematical study of localized spikes includes the
    construction of point-condensation solutions by Ni and
    Takagi \cite{NT95} and analysis of stability and
    metastability for the Gierer--Meinhardt and Gray--Scott
    systems by Wei \cite{Wei99GM,Wei99GS}, Iron, Ward and Wei \cite{IW00} and Wei and Zhang \cite{WZ01}. For the
    one-dimensional Gray--Scott model, Doelman, Gardner,
    and Kaper \cite{Doelman98} used matched asymptotic
    expansions to derive and analyze a nonlocal eigenvalue
    problem. Their subsequent work \cite{Doelman02} used
    the Evans function and a stability index to study
    spectral stability. 
    
    Motivated by the stability analysis discussed above, the following
    one-dimensional nonlocal eigenvalue problem (NLEP) on
    \(\mathbb R\) is considered:
    \begin{equation}\label{nlep}
    	L_\gamma\phi:=L_0\phi-\gamma\frac{p-1}{\int_\R w^r\,dx}w^p
    	\int_\R w^{r-1}\phi\,dx=\lambda\phi,
    	\qquad L_0:=\partial_x^2-1+pw^{p-1},
    \end{equation}
    where $p,r>1$ and $\Re\gamma\geq 0$.  Here $L_\gamma:
    H^2_{\mathrm{even}}(\mathbb R;\mathbb C)
    \rightarrow
    L^2_{\mathrm{even}}(\mathbb R;\mathbb C)$ with  
    $L^2_{\mathrm{even}}(\mathbb R;\mathbb C)
    :=\{\phi\in L^2(\mathbb R;\mathbb C):
    \phi(-x)=\phi(x)\ \text{a.e.}\}$
    and $H^2_{\mathrm{even}}(\mathbb R;\mathbb C)
    :=H^2(\mathbb R;\mathbb C)
    \cap L^2_{\mathrm{even}}(\mathbb R;\mathbb C).$   Here \(w\) is the unique positive even solution of
    \begin{align}\label{eq-w}
    w''-w+w^p=0 \quad\text{in }\mathbb R,
    \qquad w(x)\to0 \quad\text{as }|x|\to\infty.
    \end{align}
    Explicitly, 
    \begin{align}\label{explicit-profile-w}
    w(x)=
    \left(\frac{p+1}{2}\right)^{1/(p-1)}
    \sech^{2/(p-1)}
    \left(\frac{p-1}{2}x\right).
    \end{align}
    These NLEPs couple a local differential operator to weighted
    spatial integrals of the eigenfunction. The resulting operator
    is generally non-self-adjoint (unless $r=p+1$), and its spectrum need not
    be real. Establishing spectral stability therefore requires
    control of both real and nonreal eigenvalues.  The spectral properties of nonlocal perturbations have been
    studied in several settings.  Freitas \cite{Freitas94} considered a Sturm--Liouville
    operator with a bounded rank-one perturbation.
    More specifically, he established conditions ensuring
    that the spectrum is real and proved a comparison theorem.
    His later work \cite{Freitas95} gave bounds on the
    dimension of the unstable manifold and proved the existence
    of a Hopf bifurcation for a scalar nonlocal
    reaction--diffusion equation. Bose and Kriegsmann
    \cite{BK98} proved eigenvalue interlacing and an oscillation
    theorem for another class of nonlocal problems. Anastasio,
    Barreiro, and Bronski \cite{ABB17} developed a geometric
    approach to low-rank perturbations, using imaginary-axis
    eigenvalues and multiple real eigenvalues to describe
    stability regions. For an overview of nonlocal reaction--diffusion equations,
    see \cite{Freitas99}.
    
    For the radially symmetric counterpart of \eqref{nlep}
    on \(\mathbb R^N\), \(N\ge1\), the stability theory is
    well understood in certain parameter regimes, notably
    \(r=2\) and \(r=p+1\). Additional structure simplifies
    the analysis in these cases: self-adjointness for
    \(r=p+1\) and special identities permitting quadratic-form
    estimates for \(r=2\) \cite{WZ01,WWbook}.  Wei \cite{Wei99GM} obtained spectral estimates for
    the NLEP associated with single interior spikes in the
    Gierer--Meinhardt shadow system. In subsequent work,
    he established sufficient stability and instability
    criteria \cite{Wei00}. Using the basic cases \(r=2\)
    and \(r=p+1\) as starting points, Wei and Zhang
    \cite{WZ01} developed a continuation argument in \(r\)
    based on the positivity of an associated quadratic form.
    Their analysis yields explicit sufficient stability
    conditions for \(2<r<p+1\) in terms of \(p,r,\gamma\)
    and the dimension \(N\).  In one dimension, the ground state is known explicitly,
    so the spectrum of \eqref{nlep} and related problems
    can be analyzed more precisely.  For example,
    de Groen and Karadzhov \cite{DGK02} combine spectral
    decompositions with rank one resolvent formulas to
    derive a scalar characteristic equation for a nonlocal
    problem arising from the Gierer--Meinhardt shadow system.
    Their analysis treats \(r=(p+3)/2\), provides perturbative
    extensions to nearby parameter values, and yields
    instability criteria in certain regions.  For \(r=p\), Wei and Winter \cite{WW03} establish
    stability for every \(p>1\) when
    \(1<\gamma\le p/(p-1)\), with extensions to larger
    values of \(\gamma\) discussed there.
    
    Related work has considered NLEPs in which the
    nonlocal coefficient depends on the spectral parameter.
    Kolokolnikov, Sun, Ward, and Wei \cite{KSWW06}
    studied stripe instabilities and double-zero
    eigenvalues. Nec and Ward \cite{NW13} identified
    an explicitly solvable family satisfying
    \(p=2r-3\), \(r>2\), and used scalar characteristic
    equations and winding numbers to determine stability.

    We mention that the available results concern particular
    exponent relations or parameter regimes and do not cover the full parameter
    range of 
    \eqref{nlep}.  Motivated by this, we establish in this paper a complete
    classification of spectral stability for \eqref{nlep}
    on the even subspace for all \(p,r>1\) and \(\Re\gamma\ge0\).
    We determine the spectrum on the imaginary axis and count
    the eigenvalues in the right half-plane with
    algebraic multiplicity.
    
   The rank-one structure reduces the eigenvalue problem
to a scalar characteristic equation.
The main difficulty is to analyze the curve
$\gamma=1/G(i\omega)$ when the corresponding spectral
weights have mixed signs. We address this difficulty
by proving a cumulative spectral weight estimate,
which allows us to determine its geometry and count
the right-half-plane eigenvalues in each parameter region.
    
    \subsection{Main results}
    We denote the spectrum and resolvent set of
    $L_\gamma$ by $\sigma(L_\gamma)$ and
    $\rho(L_\gamma)=\mathbb C\setminus\sigma(L_\gamma)$
    respectively. Then we give the definition of spectral stability
    used throughout the paper.
    
    \begin{definition}\label{def1-lgamma}
    	We say that \(L_\gamma\) is \emph{strictly spectrally stable}
    	on the even subspace if there exists \(\delta>0\) such that
    	\[
    	\sigma(L_\gamma)
    	\subset
    	\{\lambda\in\mathbb C:\Re\lambda\le -\delta\}.
    	\]
    	Here \(\delta\) may depend on \(p\), \(r\), and \(\gamma\).
    \end{definition}
    Since $w$ is even and we can decompose $ u= u_{even}+ u_{odd}$, the odd case is reduced to local operator. So from now on we can restrict to even functions only.  For \(\lambda\in\rho(L_0)\), define
    \begin{equation}\label{def-G}
    	G(\lambda):=
    	\frac{p-1}{\int_{\mathbb R} w^r\,dx}
    	\left\langle w^{r-1},(L_0-\lambda)^{-1}w^p\right\rangle.
    \end{equation}
    Let \(N_+(L_\gamma)\) denote the number of eigenvalues of
    \(L_\gamma\) in the open right half-plane \(\Re\lambda>0\),
    counted with algebraic multiplicity.
    Our main result is stated as follows:
    \begin{theorem}\label{main}
    	For every \(p>1\) and \(r>1\), the following conclusions hold.
    	\begin{enumerate}
        \item[(i).] If \(r>(p+3)/2\), then
    \[
    N_+(L_\gamma)=
    \begin{cases}
    0, & \Re\gamma\ge 1,\\[2pt]
    0, & \displaystyle
         \exists\,\omega>0:\quad
         \Re\gamma=\Re\frac{1}{G(i\omega)},\quad
         |\Im\gamma|\ge
         \left|\Im\frac{1}{G(i\omega)}\right|,\\[4pt]
    1, & \text{otherwise}.
    \end{cases}
    \]
    $L_{\gamma}$ is strictly spectral stable  when either $\Re\gamma\ge 1$ with $\gamma\ne 1$
    or there exists \(\omega>0\) such that
    \[
    \Re\gamma=\Re\frac{1}{G(i\omega)},
    \qquad
    |\Im\gamma|>\left|\Im\frac{1}{G(i\omega)}\right|.
    \]
    Whenever the equation
    \(\Re\gamma=\Re(1/G(i\omega))\) has a solution \(\omega>0\),
    that solution is unique.
    \item[(ii).] If \(r=(p+3)/2\), then
    \[
    G(\lambda)=\frac{s_0}{s_0-\lambda},
    \qquad
    s_0=\frac{(p-1)(p+3)}{4},
    \qquad \lambda\in\rho(L_0),
    \]
    and
    \[
    N_+(L_\gamma)=
    \begin{cases}
    1, & \Re\gamma<1,\\
    0, & \Re\gamma\ge1.
    \end{cases}
    \]
    $L_{\gamma}$ is strictly spectral stable  if and only if \(\Re\gamma>1\).
    For \(\Re\gamma=1\), the only eigenvalue of \(L_\gamma\) in the
    closed right half-plane is \(s_0(1-\gamma)\), which is
    algebraically simple.
    \item[(iii).] If \(r<(p+3)/2\), then \(\Re G(i\omega)\)
    decreases strictly from \(1\) to \(0\) as \(\omega\) increases
    from \(0\) to \(\infty\). For \(\gamma\ne0\) with
    \(0<\Re(1/\gamma)<1\), let \(\omega>0\) be the unique solution of $\Re G(i\omega)=\Re(1/\gamma).$
    Then
    \[
    N_+(L_\gamma)=
    \begin{cases}
    0, & |\Im(1/\gamma)|\le \Im G(i\omega),\\
    2, & |\Im(1/\gamma)|<-\Im G(i\omega),\\
    1, & \text{otherwise}.
    \end{cases}
    \]
    For such \(\gamma\), \(L_\gamma\) is strictly spectrally stable
    if and only if $|\Im(1/\gamma)|<\Im G(i\omega).$
    If \(\gamma\notin\{0,1\}\) and \(\Re(1/\gamma)\notin(0,1)\),
    then \(N_+(L_\gamma)=1\).  
    
    Moreover, if \(r\ge(p-1)/2\), then \(\Im G(i\omega)>0\) for every
    \(\omega>0\). If \(1<r<(p-1)/2\), there is a unique \(\omega_H>0\)
    such that \(\Im G(i\omega_H)=0\), and
    \[
    \Im G(i\omega)<0 \quad (0<\omega<\omega_H),
    \qquad
    \Im G(i\omega)>0 \quad (\omega>\omega_H).
    \]
    Since \(0<\Re G(i\omega_H)<1\), we define $\gamma_H:=\frac{1}{\Re G(i\omega_H)}>1.$
    Then \(L_\gamma\) is strictly spectral stable if and only if
    \(\gamma\ne0\) and
    \[
    0<\Re(1/\gamma)<\frac{1}{\gamma_H},
    \qquad
    |\Im(1/\gamma)|<\Im G(i\omega),
    \]
    where \(\omega>0\) is the unique solution of
    \(\Re G(i\omega)=\Re(1/\gamma)\).

    		\item[(iv).]  At \(\gamma=1\), $\ker L_1=\operatorname{span}\{w\}.$  Moreover,
    		the zero eigenvalue of \(L_1\) is algebraically simple when
    		\(r\ne(p-1)/2\),
    		and has algebraic multiplicity exactly two when \(r=(p-1)/2\).  In addition,
    		\[
    		\sigma(L_1)\cap i\mathbb R=\{0\},
    		\qquad
    		N_+(L_1)=
    		\begin{cases}
    			0,& r\ge(p-1)/2,\\
    			1,& r<(p-1)/2.
    		\end{cases}
    		\]
    		
    		\item[(v).] On the critical line \(r=(p-1)/2>1\),
    		one has \(G''(0)>0\) and \(G'''(0)<0\).  As \(\gamma\rightarrow 1\) in $\C$, the two eigenvalues of \(L_\gamma\)
    		bifurcating from zero admit the following asymptotic expansions:    \begin{equation}\label{criticalexpansion}
    			\lambda_\pm=
    			\pm i\sqrt{\frac{2(\gamma-1)}{G''(0)}}
    			+\frac{G'''(0)}{3(G''(0))^2}(\gamma-1)
    			+O(|\gamma-1|^{3/2}).
    		\end{equation}
    	\end{enumerate}
    \end{theorem}
    \begin{remark}
    	Previous results \cite{WZ01,DGK02,WW03} establish stability
    	for specific exponent relations or under additional
    	parameter restrictions. For the one-dimensional problem
    	\eqref{nlep}, Theorem~\ref{main} provides necessary and
    	sufficient conditions for strict spectral stability for
    	all independent exponents \(p>1\), \(r>1\), and all $\gamma\in\C$ with
    	\(\Re\gamma\geq 0\). In particular, it covers the entire ranges
    	\(1<r<2\) and \(r>p+1\), as well as the parameter values
    	in \(2<r<p+1\) not covered by the earlier sufficient
    	conditions.
    \end{remark}

\begin{figure}[htbp]
    \centering
    \begin{subfigure}[t]{0.49\linewidth}
        \vspace{0pt}
        \centering
        \includegraphics[width=8.5cm, height=4.7cm]{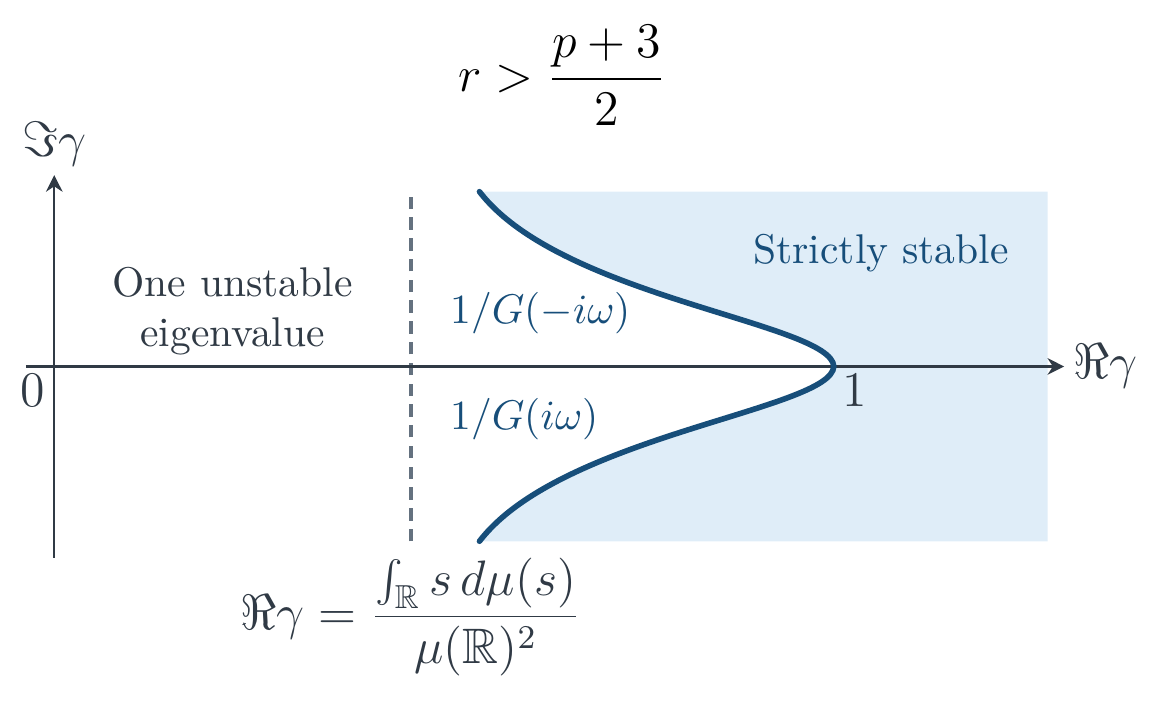}
        \caption{Part (i) for Theorem \ref{main}}
        \label{fig:stability-a}
    \end{subfigure}\hfill
    \begin{subfigure}[t]{0.49\linewidth}
        \vspace{0pt}
        \centering
        \includegraphics[width=8.5cm, height=4.7cm]{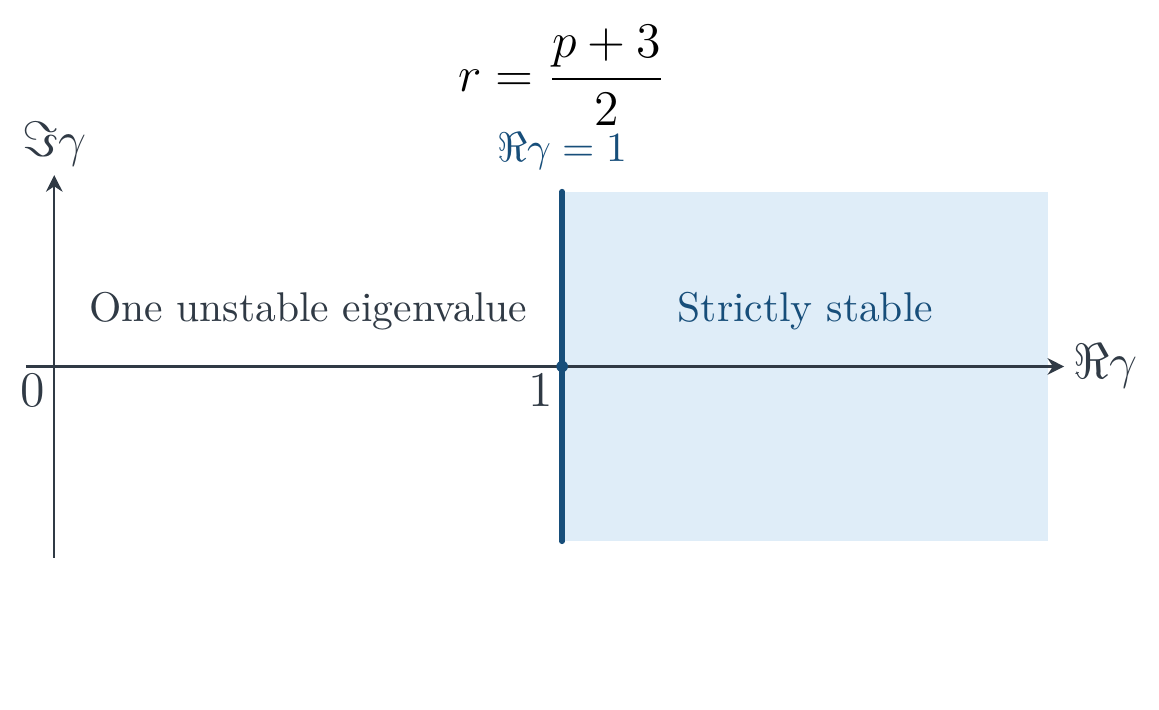}
        \caption{Part (ii) for Theorem \ref{main}}
        \label{fig:stability-b}
    \end{subfigure}

    \medskip

    \begin{subfigure}[t]{0.49\linewidth}
        \vspace{0pt}
        \centering
        \includegraphics[width=8.5cm, height=4.7cm]{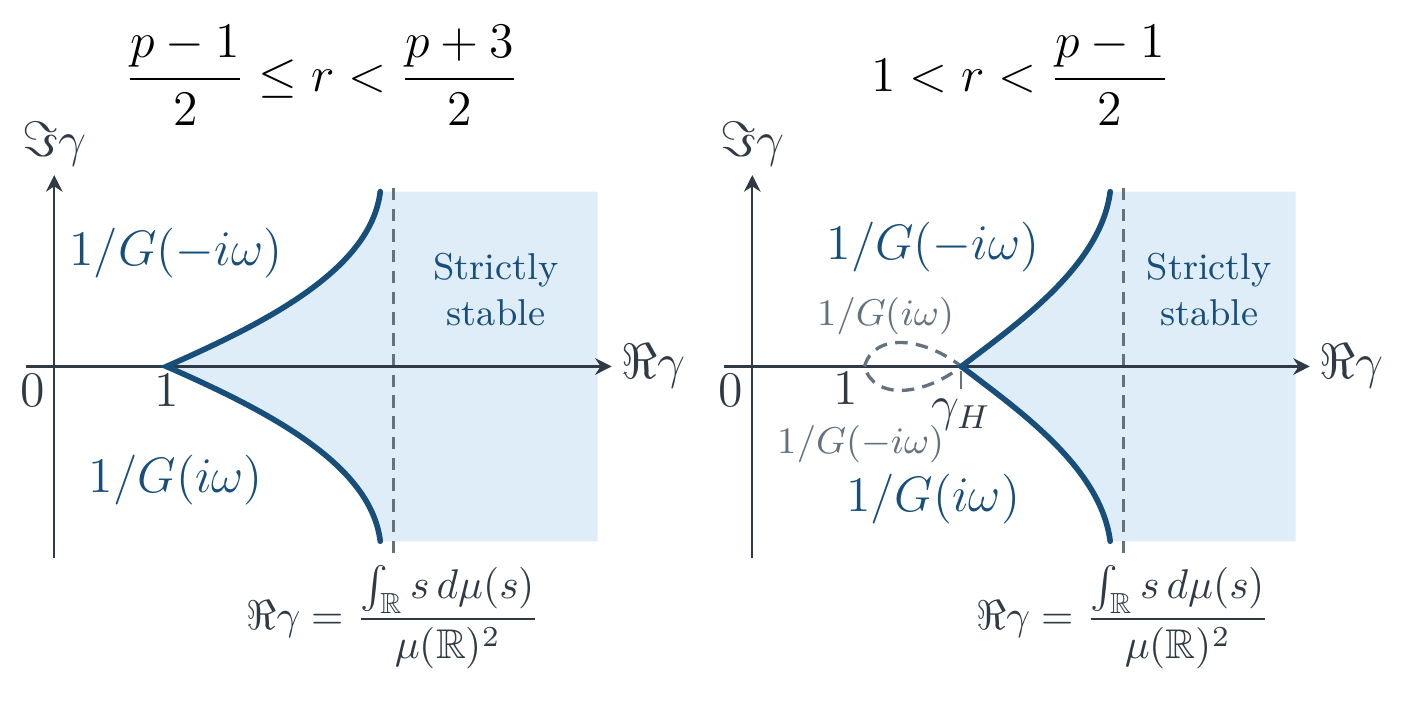}
        \caption{Part (iii) for Theorem \ref{main} in the complex \(\gamma\)-plane}
        \label{fig:stability-c}
    \end{subfigure}\hfill
    \begin{subfigure}[t]{0.49\linewidth}
        \vspace{0pt}
        \centering
        \includegraphics[width=8.5cm, height=4.7cm]{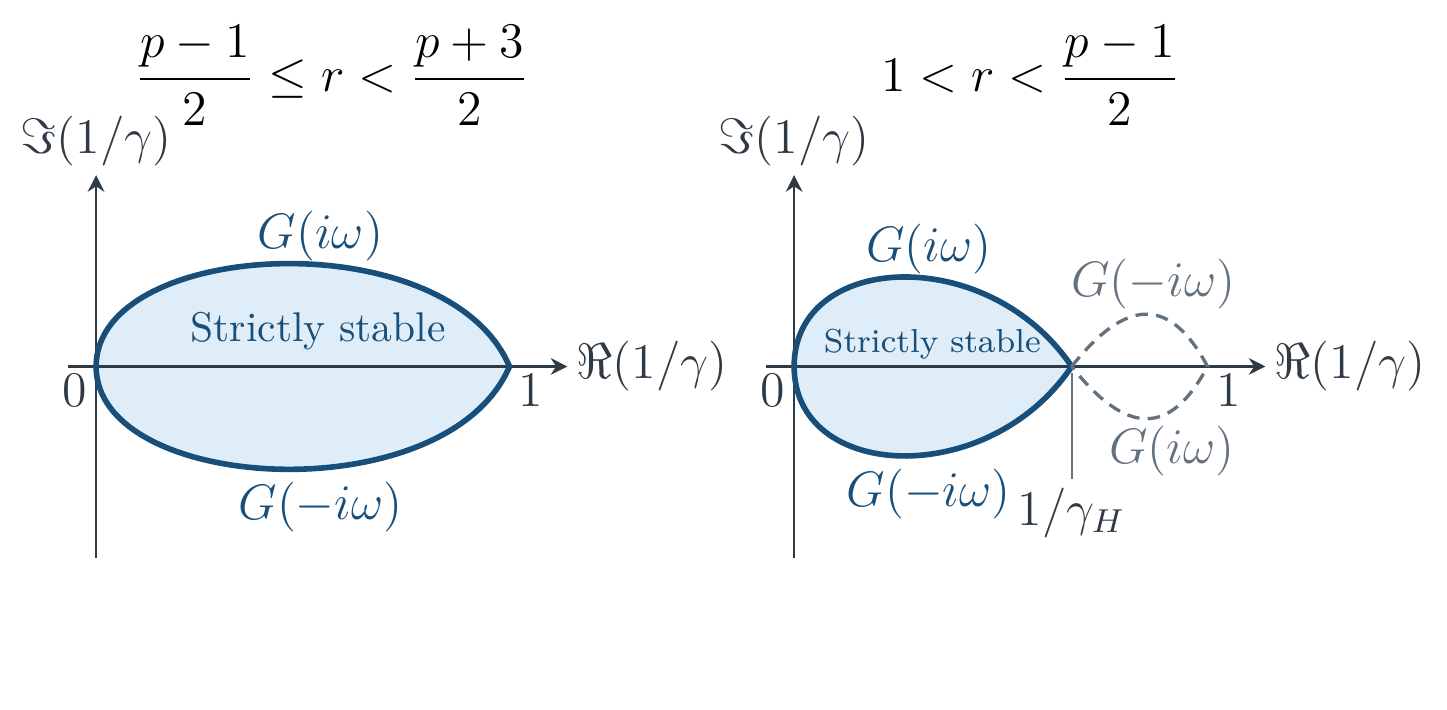}
        \caption{Part (iii) for Theorem \ref{main} in the complex
    \(1/\gamma\)-plane}
        \label{fig:stability-d}
    \end{subfigure}

    \caption{Schematic diagram: spectral stability regions of $L_{\gamma}$.}
    \label{fig:stability}
\end{figure}

Theorem~\ref{main} exhibits necessary and sufficient conditions for
the strict spectral stability of \(L_\gamma\) in different parameter
regimes, with strict spectral stability as in
Definition~\ref{def1-lgamma}. In Figure~\ref{fig:stability}(a), where \(r>(p+3)/2\), the stability
boundary consists of two conjugate arcs that open to the left from
their common endpoint \(\gamma=1\). The strictly stable
region lies to the right of the boundary. Whereas, in
Figure~\ref{fig:stability}(c), where \(r<(p+3)/2\), the two conjugate
arcs open to the right, and the strictly stable region lies between
them. Their common endpoint is \(\gamma=1\) when
\((p-1)/2\le r<(p+3)/2\) and \(\gamma=\gamma_H>1\) when
\(1<r<(p-1)/2\).
{In both panels, the two unbounded arcs approach the
vertical asymptote
\(\Re\gamma=\frac{\int_{\mathbb R}s\,d\mu(s)}{\mu(\mathbb R)^2}\)
as \(\omega\to\infty\).} In particular, when $1<r<(p-1)/2$, the dashed arcs
$1/G(\pm i\omega)$, $0\le\omega\le\omega_H$,
connect $\gamma=1$ to $\gamma=\gamma_H$
and enclose a bounded region where
$N_+(L_\gamma)=2$.  As illustrated in Figure~\ref{fig:stability}(b), at \(r=(p+3)/2\)
the vertical line \(\Re\gamma=1\) separates the strictly stable
half-plane \(\Re\gamma>1\) from the unstable half-plane
\(\Re\gamma<1\).  We remark that  Figures~\ref{fig:stability}(c) and~\ref{fig:stability}(d) show the same
stability classification in the \(\gamma\)- and \(1/\gamma\)-planes,
respectively. The reciprocal coordinate simplifies the proof because
the curves \(G(\pm i\omega)\), \(\omega>0\), are graphs over the real axis.

Restricting $\gamma$ to real values in Theorem~\ref{main}, we obtain
 the following corollary.
\begin{corollary}\label{cor:real-gamma}
Let \(p,r>1\) and \(\gamma\ge0\) be real. For \(0\le\gamma<1\),
\[
N_+(L_\gamma)=1,\qquad
\sigma(L_\gamma)\cap i\mathbb R=\varnothing.
\]
If \(r\ge(p-1)/2\), then \(L_\gamma\) is strictly spectral
stable if and only if \(\gamma>1\).

If \(1<r<(p-1)/2\), let \(\omega_H>0\) be the unique zero
of \(\Im G(i\omega)\) for \(\omega>0\), and set
\(\gamma_H=1/\Re G(i\omega_H)>1\). Then
\[
N_+(L_\gamma)=
\begin{cases}
2,&1<\gamma<\gamma_H,\\
0,&\gamma\ge\gamma_H.
\end{cases}
\]
Strict spectral stability holds if and only if \(\gamma>\gamma_H\).
At \(\gamma=\gamma_H\), the only spectrum in the closed right
half-plane is the algebraically simple pair \(\pm i\omega_H\),
which crosses the imaginary axis transversely from right to left
as \(\gamma\) increases.
\end{corollary}

        \section{Preliminaries:  Spectrum of the Local Operator}
    In this section, we first recall some important results on the spectrum of the local operator \(L_0\) given by \(L_0:=\partial_x^2-1+pw^{p-1}\).  Set
    \begin{equation}\label{scalednotation}
    	\ell=\frac{p+1}{p-1}>1,\qquad k=\frac1{\ell-1},\qquad
    	q=\frac{2(r-1)}{p-1}>0,\qquad f_u(z)=\sech^u z.
    \end{equation}
    With $z=kx$,  $L_0$ becomes
    \begin{equation}\label{scaledoperator}
    	L_0=k^2\bigl(\partial_z^2+\ell(\ell+1)\sech^2 z\bigr)-1,
    \end{equation} and the function \(G\) defined in \eqref{def-G} can be written as
    \begin{equation}\label{scalingconstant}
    	G(\lambda)=c\langle{f_q},{(  L_0 -\lambda)^{-1}f_{\ell+1}}\rangle,\qquad c=\frac{p^2-1}{2\int_0^\infty f_{\frac{2r}{p-1}}(z)\,dz}.
    \end{equation}
    The spectral properties of the local operator $L_0$  are presented as follows.
    \begin{lemma}\label{local}
    	For \(m\in\mathbb C\), the solution of
    	\begin{align}\label{crucial-ODEs}
    		\psi''+\ell(\ell+1)\sech^2 z\,\psi=m^2\psi,
    		\qquad z>0,
    	\end{align}
    	satisfying \(\psi_m'(0)=0\) and normalized by \(\psi_m(0)=1\)
    	is given by
    	\begin{equation}\label{regularsolution}
    		\psi_m(z)=f_m(z)\,
    		{}_2F_1\left(
    		\frac{m-\ell}{2},\frac{m+\ell+1}{2};
    		\frac12;\tanh^2 z
    		\right),
    	\end{equation}
    	where ${}_2F_1$ is the Gauss hypergeometric function.
    	
    	Consider the eigenvalue problem for \(L_0\):
    	\begin{align}\label{llep-s-lemma}
    		L_0\phi=s\phi,\qquad
    		0\not\equiv\phi\in L^2(\mathbb R;\mathbb C),\qquad s\in\mathbb R.
    	\end{align}
    	Set \(J=\lceil\ell/2\rceil-1\), the eigenvalues on the even subspace are given by
    	\begin{align}\label{sj-k-m-l-j}
    		s_j=k^2m_j^2-1,\qquad
    		m_j=\ell-2j>0,\qquad j=0,\ldots,J,
    	\end{align}
    	with corresponding eigenfunctions \(\psi_{m_j}(kx)\in L^2\).
    	In particular, \(s_0>0\), with positive eigenfunction
    	\(f_\ell(kx)\). For $j=1,\ldots,J$, we have $s_j<0$.
    	The continuous spectrum of \(L_0\) on the even subspace
    	is \((-\infty,-1]\), and \(0\in\rho(L_0)\).
    \end{lemma}
    
    \begin{proof}
    	Formula \eqref{regularsolution} follows from the standard even
    	solution of the modified P\"oschl--Teller equation, see 
    	(39.10a) in \cite{Flugge99}.  Thanks to (39.25) in 
    	\cite{Flugge99}, one has  
    	\begin{align}\label{mj-ell-2j}
    	m_j=\ell-2j>0,\qquad
    	j=0,\ldots,\lceil\ell/2\rceil-1.
    	\end{align}
    	Using $s=k^2m^2-1$, we obtain the eigenvalues $s_j$ and
    the corresponding eigenfunctions $\psi_{m_j}(kx)$
    in \eqref{sj-k-m-l-j}.

    For $j=0$, \eqref{regularsolution} gives
    $\psi_\ell=f_\ell>0$. The corresponding eigenvalue is
    \begin{align}\label{s0-beta-p-3}
        s_0=k^2\ell^2-1
        =\frac{(p-1)(p+3)}{4}>0.
    \end{align}  For \(j\ge1\), since $0<m_j\le\ell-2,$ we have \(s_j<0\).
    
    Finally, by using  Theorem~9.38 in \cite{teschl}, the continuous spectrum of $-\partial_z^2-\ell(\ell+1)\sech^2 z$ is $[0,\infty)$. Hence, the continuous spectrum of $L_0$ on the even subspace is $(-\infty,-1]$. Since radial functions in one dimension are even, Lemma~4(1) of \cite{WZ01} implies that $0\in\rho(L_0)$.
    \end{proof}
    We now introduce the spectral transform associated with \(L_0\) for $f_u$.   For \(j=0,\ldots,J\) and \(u>0\), we define the discrete
    transform of \(f_u\) by
    \begin{equation}\label{boundtransform}
    	J_j(u):=\int_0^\infty f_u\psi_{m_j}\,dz.
    \end{equation}
    Similarly, for each \(u>0\), the continuous transform of \(f_u\) is defined by
    \begin{equation}\label{continuumtransform}
    	T_u(\xi):=\int_0^\infty f_u(z)\psi_{i\xi}(z)\,dz,
    	\qquad \xi>0.
    \end{equation} 
    
    
    The following lemma gives explicit formulas for the discrete
    and continuous spectral weights  of
    \(G\) and determines their signs.
    \begin{lemma}\label{spectralweights}
    	Let \(E\) be the projection-valued spectral measure of \(L_0\). For each Borel set \(B\subset\mathbb R\), define
    	\begin{equation}\label{crossmeasure}
    		\mu(B):=c\left\langle f_q,E(B)f_{\ell+1}\right\rangle.
    	\end{equation}
    	Then \(\mu\) is a finite signed measure with
    	\begin{align}\label{tv-bound-G-spectral}
    		\|\mu\|_{\mathrm{TV}}
    		\le c\|f_q\|_{L^2(0,\infty)}
    		\|f_{\ell+1}\|_{L^2(0,\infty)},~~~G(\lambda)=\int_{\mathbb R}\frac{d\mu(s)}{s-\lambda},
    	\end{align}
    	for  $\lambda\in\rho(L_0)$ and $c$ is given in \eqref{scalingconstant}.  For \(j=0,\ldots,J\), the spectral weight of \(\mu\) at \(s_j\) is
    	\begin{equation}\label{atommass}
    		\mu_j:=\mu(\{s_j\})
    		=c\frac{J_j(q)J_j(\ell+1)}{N_j},~N_j:=\int_0^\infty|\psi_{m_j}|^2\,dz
    		=\frac{\pi\Gamma(j+1)\Gamma(\ell-j+1)}
    		{2(\ell-2j)\Gamma(j+1/2)\Gamma(\ell-j+1/2)}
    	\end{equation}
    	with $J_j$  satisfying
    	\begin{equation}\label{boundtransform-explicit}
    		J_j(u)
    		=\frac{\sqrt\pi}{2}
    		\frac{\Gamma((u+\ell)/2-j)}{\Gamma((u+\ell+1)/2)}
    		\left(\frac{u-\ell}{2}\right)_j
    	\end{equation}
    	and  \((a)_j\) denotes the Pochhammer symbol.
    	The continuous part  \(\mu_{\mathrm c}\) of \(\mu\) is given by
        \begin{align}\label{mu-c-B}
    	\mu_{\mathrm c}(B)
    	=c\int_0^\infty
    	\mathbf{1}_B(-1-k^2\xi^2)
    	T_q(\xi)T_{\ell+1}(\xi)\,d\rho_{\mathrm c}(\xi),
    	\end{align}
    	where the continuous spectral measure corresponding to
    	\(\psi_{i\xi}(0)=1\) and \(\psi_{i\xi}'(0)=0\) is
    	\begin{equation}\label{weylmeasure}
    		d\rho_{\mathrm c}(\xi)
    		=
    		\frac{
    			\left|
    			\Gamma\!\left(\frac{-\ell+i\xi}{2}\right)
    			\Gamma\!\left(\frac{\ell+1+i\xi}{2}\right)
    			\right|^2
    		}{
    			2\pi^2|\Gamma(i\xi)|^2
    		}\,d\xi,
    		\qquad \xi>0.
    	\end{equation} 
    	The continuous transform \eqref{continuumtransform} is given explicitly by \begin{equation}\label{continuumtransform-explicit}
    		T_u(\xi)
    		=\frac{\sqrt\pi}{2}
    		\frac{|\Gamma((u+i\xi)/2)|^2}
    		{\Gamma((u+\ell+1)/2)\Gamma((u-\ell)/2)},
    		\qquad \xi>0.
    	\end{equation}
    	In particular, \(\mu_0>0\),
    	\(J_j(\ell+1)>0\), and \(T_{\ell+1}(\xi)>0\). Consequently,
    	\begin{equation}\label{signs}
    		\operatorname{sgn}\mu_j
    		=\operatorname{sgn}\left(\frac{q-\ell}{2}\right)_j,
    		\qquad
    		\operatorname{sgn}T_q(\xi)
    		=\operatorname{sgn}\frac{1}{\Gamma((q-\ell)/2)},
    		\qquad \xi>0.
    	\end{equation}
    	
    \end{lemma}
    
    \begin{proof}
    The operator $L_0$ is self-adjoint on
    $L^2_{\mathrm{even}}(\mathbb R;\mathbb C)$. By the spectral theorem (Theorem~3.7 in \cite{teschl}),
    $\mu$  is a finite
    complex measure. Since $L_0$ commutes with complex
    conjugation and $f_q,f_{\ell+1}$ are real-valued,
    $\mu$ is real-valued.
    The total variation bound in \eqref{tv-bound-G-spectral}
    follows from Cauchy--Schwarz and the orthogonality
    of the spectral projections.
    The spectral theorem also gives
    \begin{align}\label{spectral-representation-G}
    	G(\lambda)
    	&=\int_{\mathbb R}\frac{d\mu(s)}{s-\lambda},\quad
    \lambda\in\rho(L_0),\nonumber\\
    L_0f&=\int_{\mathbb R}s\,dE(s)f,
    \quad f\in H^2_{\mathrm{even}}(\mathbb R;\C),
    \end{align}
    
    By Sturm--Liouville theory, its eigenvalues on
    the even subspace are simple. Consequently,
    	\[
    	E(\{s_j\})f
    	=\frac{\langle\psi_{m_j},f\rangle}{N_j}\,\psi_{m_j},~N_j=\int_0^\infty|\psi_{m_j}(z)|^2\,dz.
    	\]
    	and hence
    	\[
    	\mu_j
    	=c\langle f_q,E(\{s_j\})f_{\ell+1}\rangle
    	=c\frac{J_j(q)J_j(\ell+1)}{N_j}.
    	\]
    The formula \eqref{boundtransform-explicit} follows
    from \eqref{boundtransform} and \eqref{regularsolution}.

    For $m=m_j=\ell-2j$, the Jacobi polynomial
    representation of \eqref{regularsolution} and the
    substitution $t=\tanh^2z$ give
    	\[
    	N_j=\frac12\left(\frac{j!}{(1/2)_j}\right)^2
    	\int_0^1 t^{-1/2}(1-t)^{m-1}
    	\left[P_j^{(-1/2,m)}(1-2t)\right]^2\,dt,
    	\]
        where \(P_j^{(-1/2,m)}\) is the Jacobi polynomial of degree \(j\)
    with parameters \(-1/2\) and \(m\).
    	Applying the Jacobi integral formula
    	((2.21) in \cite{Nieto78}), one has
    	\begin{align}\label{Ni-normalization}
    		N_j=
    		\frac{\pi\Gamma(j+1)\Gamma(\ell-j+1)}
    		{2(\ell-2j)\Gamma(j+1/2)\Gamma(\ell-j+1/2)}.
    	\end{align}

    	For \(\xi>0\), the solution to \eqref{crucial-ODEs} has the asymptotic
    	expansion
    	\[
    	\psi_{i\xi}(z)
    	=A(i\xi)e^{i\xi z}
    	+\overline{A(i\xi)}e^{-i\xi z}+o(1),
    	\qquad z\to\infty,
    	\]
        where $A$ is a complex constant.  Thanks to Weyl--Titchmarsh theory
    	(Theorem~9.17 and (9.63) in \cite{teschl}), one finds 
    	\[
    	d\rho_{\mathrm c}(\xi)
    	=\frac{d\xi}{2\pi|A(i\xi)|^2},
    	\]
    	Here \(\rho_{\mathrm c}\) is the continuous spectral measure
    	of  \(L_0\).  It then follows from  \(\psi_{i\xi}(0)=1\) and
    	\(\psi_{i\xi}'(0)=0\) that \eqref{weylmeasure} holds.  By the spectral representation of $L_0$ and the definition of $\mu$, we obtain
    	\[
    	\mu_{\mathrm c}(B)
    	=c\int_0^\infty
    	\mathbf1_B(-1-k^2\xi^2)
    	T_q(\xi)T_{\ell+1}(\xi)\,d\rho_{\mathrm c}(\xi).
    	\]
    Substituting \eqref{regularsolution} into
    \eqref{continuumtransform}, we obtain
    \eqref{continuumtransform-explicit}.  Moreover,
    	\[
    	\begin{aligned}
    		&c\int_0^\infty
    		\left|\mathbf1_B(-1-k^2\xi^2)
    		T_q(\xi)T_{\ell+1}(\xi)\right|\,d\rho_{\mathrm c}(\xi)\\
    		&\quad\le
    		c\left(\int_0^\infty|T_q(\xi)|^2\,d\rho_{\mathrm c}(\xi)\right)^{1/2}
    		\left(\int_0^\infty|T_{\ell+1}(\xi)|^2\,d\rho_{\mathrm c}(\xi)\right)^{1/2}\\
    		&\quad\le
    		c\|f_q\|_{L^2(0,\infty)}
    		\|f_{\ell+1}\|_{L^2(0,\infty)}
    		<\infty,
    	\end{aligned}
    	\]
    	where we used \eqref{continuumtransform-explicit} and \eqref{boundtransform-explicit}.
    	
    Finally, by using \eqref{boundtransform-explicit}, we get
    $J_j(\ell+1)>0$ and the sign of $J_j(q)$.
    Since $c>0$ and $N_j>0$, \eqref{atommass}
    yields the sign formula for $\mu_j$ in
    \eqref{signs}. In particular, $\mu_0>0$.
    Similarly, thanks to \eqref{continuumtransform-explicit}, one obtains
    $T_{\ell+1}(\xi)>0$ and the sign formula for
    $T_q(\xi)$ in \eqref{signs}.
    
    \end{proof}


    \section{Rank one Reduction and Moment Identities}
    In this section, we reduce \eqref{nlep} to a scalar
    characteristic equation and relate the orders of its
    zeros to the algebraic multiplicities of the corresponding
    eigenvalues. We also establish a formula for counting
    eigenvalues in the right half-plane and derive moment
    identities for \(G\).  Define the perturbation
    \begin{align}\label{rank-1-formula}
    R\phi:=\frac{p-1}{\int_\R w^r\,dx}
    w^p\langle w^{r-1},\phi\rangle
    \end{align}
    and rewrite $L_\gamma=L_0-\gamma R$, where $L_{\gamma}$ is given by \eqref{nlep}.
    
    \begin{lemma}\label{rankone}
    The operator $R$  is bounded and has rank one.  
    For $\lambda\in\rho(L_0)$, $\lambda$ is an eigenvalue
    of $L_\gamma$ if and only if
    $d_\gamma(\lambda)=1-\gamma G(\lambda)=0$.
    The algebraic multiplicity of $\lambda$ equals
    its multiplicity as a zero of $d_\gamma$.   In particular, for every  $\gamma\not=0$, $s_0\in\rho(L_{\gamma})$.
    
    For every compact set $K\subset\{\gamma\in\mathbb C:\Re\gamma\ge0\}$,
    there exists $M>s_0$ with $s_0$ given by \eqref{s0-beta-p-3} such that, for all
    $\gamma\in K$,
    \[
    \sigma(L_\gamma)\cap\{\Re\lambda\ge0\}
    \subset
    \{\lambda\in\mathbb C:
    0\le\Re\lambda<M,\ |\Im\lambda|<M\}.
    \]
    This part of the spectrum consists of finitely
    many isolated eigenvalues, each one has finite algebraic
    multiplicity.
     Furthermore, for every  $\gamma_0\in\mathbb C$ with $\Re\gamma_0\ge0$,
    	\[
    	\sigma(L_{\gamma_0})\cap i\mathbb R=\varnothing
    	\quad\Longrightarrow\quad
    	\exists\,\varepsilon>0:\quad
    	N_+(L_\gamma)=N_+(L_{\gamma_0})
    	\quad
    	\text{for all } \gamma\in\mathbb C,~\Re\gamma\ge0
    	\text{ with }|\gamma-\gamma_0|<\varepsilon.
    	\]
    	
    \end{lemma}

    \begin{proof}
    
    Since $w^p,w^{r-1}\in L^2$, the operator $R$
    is bounded and has rank one. In particular $R$ is compact,
    and
    \[
    \sigma_{\mathrm{ess}}(L_\gamma)
    =\sigma_{\mathrm{ess}}(L_0)=(-\infty,-1],
    \]
    where \(\sigma_{\mathrm{ess}}(L_{\gamma})\) denotes the Fredholm
    essential spectrum of \(L_{\gamma}\).
    Noting that \(L_\gamma-\lambda\) is invertible for sufficiently
    large positive \(\lambda\), by analytic Fredholm theory (Theorem~1.9 of Chapter~VII in \cite{Kato}), the spectrum outside
    $(-\infty,-1]$ consists of isolated eigenvalues
    of finite algebraic multiplicity.
    
    For $\lambda\in\rho(L_0)$ with
    $d_\gamma(\lambda)\ne0$, we obtain from \eqref{rank-1-formula} that 
    \begin{equation}\label{resolvent}
    	(L_\gamma-\lambda)^{-1}
    	=(L_0-\lambda)^{-1}
    	+\frac{\gamma}{d_\gamma(\lambda)}
    	(L_0-\lambda)^{-1}R(L_0-\lambda)^{-1},
    \end{equation}
    where $d_\gamma(\lambda)=1-\gamma G(\lambda)$. Since $R$ has rank one, we have
    \[
    \det\bigl(I-\gamma R(L_0-\lambda)^{-1}\bigr)
    =1-\gamma G(\lambda)=d_\gamma(\lambda),
    \qquad \lambda\in\rho(L_0).
    \]
    The Weinstein--Aronszajn formula (Theorem~6.2 of Chapter~IV  in \cite{Kato})
    then shows that
    $\lambda$ is an eigenvalue of
    $L_\gamma$ if and only if $d_\gamma(\lambda)=0$.  Moreover,
    its algebraic multiplicity equals the order
    of $d_\gamma$ at $\lambda$.

    	Since $L_0$ is self-adjoint and $R$ is bounded,
    	Theorem~IV.3.17 in \cite{Kato} implies
    	\[
    	\sigma(L_\gamma)
    	\subset
    	\left\{
    	\lambda\in\mathbb C:
    	\operatorname{dist}(\lambda,\sigma(L_0))
    	\le|\gamma|\|R\|
    	\right\}.
    	\]
    	Hence, for every compact set $K\subset\{\gamma\in\mathbb C:\Re\gamma\ge0\}$, uniformly in $\gamma\in K$,
    	\begin{align}\label{uniform enclosure}
    		\sigma(L_\gamma)\cap\{\Re\lambda\ge0\}
    		\subset
    		\left\{
    		\lambda\in\mathbb C:
    		0\le\Re\lambda\le s_0+(\max_{\zeta\in K} |\zeta|)\|R\|,
    		\quad
    		|\Im\lambda|\le(\max_{\zeta\in K} |\zeta|)\|R\|
    		\right\}.
    	\end{align}
    	Since this enclosure is disjoint from
    	$\sigma_{\mathrm{ess}}(L_\gamma)=(-\infty,-1]$, we choose $\mathcal C$ as a counterclockwise rectangular contour
    enclosing all eigenvalues with positive real part,
    with its left edge on the imaginary axis, then
    	continuity of Riesz projections $-\frac1{2\pi i}
    \oint_{\mathcal C}(L_\gamma-\lambda I)^{-1}\,d\lambda$ implies that
    	$N_+(L_\gamma)$ is locally constant wherever
    	$\sigma(L_\gamma)\cap i\mathbb R=\varnothing$.
    	
    \end{proof}
    To determine changes in \(N_+(L_\gamma)\), it suffices
    to track the eigenvalues near the imaginary axis.
    This is established in the following lemma.
    \begin{lemma}\label{localcount}
    	Fix $\gamma_*\in\mathbb C$ with $\Re\gamma_*\ge0$ and suppose $\sigma(L_{\gamma_*})\cap i\mathbb R=\{\lambda_1,\ldots,\lambda_m\},$
    	where each \(\lambda_j\) is an isolated eigenvalue of finite
    	algebraic multiplicity. Let \(D_j\) be disks centered at
    	\(\lambda_j\), with pairwise disjoint closures, such that
    	\[
    	\overline{D_j}\cap\sigma(L_{\gamma_*})=\{\lambda_j\},
    	\qquad j=1,\ldots,m.
    	\]
    	Then there exists \(\delta>0\) such that,
    	for every $\gamma\in\mathbb C$ with $\Re\gamma\ge0$ satisfying \(|\gamma-\gamma_*|<\delta\),
    	\[
    	\partial D_j\subset\rho(L_\gamma),
    	\qquad
    	n_{D_j}(L_\gamma)=n_{D_j}(L_{\gamma_*})<\infty,
    	\qquad j=1,\ldots,m,
    	\]
    Here $n_U(L)$ denotes the total algebraic
    multiplicity of the eigenvalues of $L$ in $U$. Moreover,
    	\begin{equation}\label{localcountformula}
    		N_+(L_\gamma)
    		-\sum_{j=1}^m
    		n_{D_j\cap\{\Re\lambda>0\}}(L_\gamma)
    		=N_+(L_{\gamma_*}),
    	\end{equation}
    	where \(N_+(L_\gamma)\) denote the number of eigenvalues of
    	\(L_\gamma\) in the open right half-plane \(\Re\lambda>0\).
    \end{lemma}

    \begin{proof}
    By using  Theorem 3.16 of Chapter~IV in \cite{Kato}, we obtain that the Riesz
    projections associated with $\partial D_j$
    have constant finite ranks for $\gamma$
    sufficiently close to $\gamma_*$. To separate the eigenvalues near the imaginary
    axis from the remaining right-half-plane spectrum, we
    choose $M>0$ sufficiently large thanks to
    \eqref{uniform enclosure} and set
    \[
    U=
    \{\lambda\in\mathbb C:
    0<\Re\lambda<M,\ |\Im\lambda|<M\}
    \setminus\bigcup_{j=1}^m\overline{D_j}.
    \]
    	Then $\partial U\subset\rho(L_{\gamma_*})$.
    	Applying the same theorem in \cite{Kato} to $U$, we obtain
    	\[
    	n_U(L_\gamma)
    	=n_U(L_{\gamma_*})
    	=N_+(L_{\gamma_*})
    	\]
    	for $\gamma$ sufficiently close to $\gamma_*$.
    	Since
    	\[
    	N_+(L_\gamma)
    	=n_U(L_\gamma)
    	+\sum_{j=1}^m
    	n_{D_j\cap\{\Re\lambda>0\}}(L_\gamma),
    	\]
    	this proves \eqref{localcountformula}.
    \end{proof}
    
    We next derive several moment identities.
    \begin{lemma}\label{moments}
    	For every \(p>1,r>1\),
    	\begin{equation}\label{exactmoments}
    		\begin{gathered}
    			G(0)=1,\qquad G'(0)=\frac1{p-1}-\frac1{2r},\\
    			\mu(\R)=\frac{p-1}{\int_\R w^r\,dx}\int_\R w^{p+r-1}\,dx=\frac{(p-1)(p+1)r}{2r+p-1}>0,\\
    \int_{\mathbb R}s\,d\mu(s)
    =
    \frac{(p-1)^2(p+1)r}{2r+p-1}
    \frac{p^2+2p+2r-3}{3p+2r-3}>0.
    		\end{gathered}
    	\end{equation}
    In addition, $0\in\sigma(L_\gamma)$ if and only if
    $\gamma=1$, and  $\dim\ker L_1=1$.
     Moreover, if \(r\ne(p-1)/2\), zero is an algebraically simple
    	eigenvalue of \(L_1\). For \(\gamma\) sufficiently close to
    	\(1\) in $\C$, the unique eigenvalue \(\lambda(\gamma)\) of
    	\(L_\gamma\) near zero satisfies
    	\begin{equation}\label{simplebranch}
    		\lambda(\gamma)
    		=-\frac{2r(p-1)}{2r-p+1}(\gamma-1)
    		+O\bigl(|\gamma-1|^2\bigr).
    	\end{equation}
    \end{lemma}
    
    \begin{proof}
    In light of the equation for $w$ given by \eqref{eq-w}, one has
    	\[
    	L_0w=(p-1)w^p,\qquad
    	L_0\left(\frac{w}{p-1}+\frac{xw'}2\right)=w.
    	\]
    The first identity yields $G(0)=1$. Differentiating $G$ and using the second identity, we obtain
    	\[
    	G'(0)=\frac1{\int_\R w^r}\int w^{r-1}
    	\left(\frac{w}{p-1}+\frac{xw'}2\right)
    	=\frac1{p-1}-\frac1{2r}.
    	\]
        By using  \eqref{eq-w}, we show that 
     \begin{align}\label{moment-pre-mode-1}
    (w')^2=w^2-\frac{2}{p+1}w^{p+1},~L_0w^p=(p-1)w^p+p(p-1)w^{p-2}(w')^2.
    \end{align}
    Moreover, thanks to \eqref{explicit-profile-w}, one has 
    \begin{align}\label{useful-ratio-moment}
    \frac{\int_{\mathbb R}w^{u+p-1}\,dx}
         {\int_{\mathbb R}w^u\,dx}
    =\frac{(p+1)u}{2u+p-1},
    \qquad u>0.
    \end{align}
    Invoking \eqref{spectral-representation-G} and \eqref{crossmeasure}, one finds 
    \begin{align}\label{0-1-modes}
    \mu(\mathbb R)
    =\frac{p-1}{\int_{\mathbb R}w^r\,dx}
     \int_{\mathbb R}w^{p+r-1}\,dx,~\int_{\mathbb R}s\,d\mu(s)
    =\frac{p-1}{\int_{\mathbb R}w^r\,dx}
     \left\langle w^{r-1},L_0w^p\right\rangle.
    \end{align}
    Then, substituting \eqref{moment-pre-mode-1} and \eqref{useful-ratio-moment} into \eqref{0-1-modes}, we get the explicit formulas of $\mu(\mathbb R)$ and $\int_{\mathbb R}s\,d\mu(s)$  shown  in \eqref{exactmoments}.

    In addition, 
    \[
    d_\gamma(0)=1-\gamma G(0)=1-\gamma,
    \]
    so zero is an eigenvalue exactly when \(\gamma=1\) thanks to Lemma \ref{rankone}.  Since $L_0^{-1}w^p=w/(p-1)$, the equation
    $L_1\phi=0$ is equivalent to $\phi=
    \frac{\langle w^{r-1},\phi\rangle}
    {\int_{\mathbb R}w^r\,dx}\,w.$
    This shows that $\ker L_1=\operatorname{span}\{w\}$.
    
    	If \(r\ne(p-1)/2\), then
    	\[
    	d_1'(0)=-G'(0)
    	=-\frac{2r-p+1}{2r(p-1)}\ne0.
    	\]
    	By Lemma~\ref{rankone}, we obtain $\ker(L_1^2)=\ker L_1.$
    	The analytic implicit function theorem implies  that there is a unique analytic function \(\lambda(\gamma)\) near \(\gamma=1\)
    	satisfying
    	\[
    	d_\gamma(\lambda(\gamma))=0,
    	\qquad \lambda(1)=0.
    	\]
    	Differentiating this identity at \(\gamma=1\), we obtain
    	\[
    	-G(0)-G'(0)\lambda'(1)=0,
    	\qquad
    	\lambda'(1)
    	=-\frac{1}{G'(0)}
    	=-\frac{2r(p-1)}{2r-p+1},
    	\]
    	which completes the proof of
    	\eqref{simplebranch}.
    \end{proof}

\section{Sign and Monotonicity Estimates for Eigenvalue Crossings}
 {
In this section, we establish cumulative spectral weight
estimates in the mixed-sign regime
$1<p<3$, $1<r<(5-p)/2$, or equivalently,
$\ell>2$, $0<q<\ell-2$.
By Lemmas~\ref{rankone} and~\ref{localcount},
$N_+(L_\gamma)$ can change only through
imaginary-axis crossings, characterized by
$1/\gamma=G(i\omega)$.
To find the stability regions in the
$1/\gamma$-plane, we study the monotonicity of
$\Re G(i\omega)$ and the sign of $\Im G(i\omega)$,
which determine the geometry of the conjugate
curves $G(\pm i\omega)$.
By \eqref{spectral-representation-G}, for $\omega>0$,
\begin{align}\label{def-H-y-s2-y}
\begin{aligned}
\Re G(i\omega)
&=\frac{s_0\mu_0}{s_0^2+\omega^2}
+\int_{(-\infty,0)}
\frac{s\,d\mu(s)}{s^2+\omega^2},\\
H(\omega^2):=\frac{\Im G(i\omega)}{\omega}
&=\frac{\mu_0}{s_0^2+\omega^2}
+\frac{\mu_1}{s_1^2+\omega^2}
+\int_{(-\infty,0)\setminus\{s_1\}}
\frac{d\mu(s)}{s^2+\omega^2}.
\end{aligned}
\end{align}
}
The strict monotonicity of
\(\Re G(i\omega)\) for \(\omega>0\) ensures that
\(\Re G(i\omega)=\Re(1/\gamma)\) has at most one positive solution.
At that frequency, the intersections have imaginary parts
\(\pm\Im G(i\omega)\), and the sign of \(\Im G(i\omega)\) determines
which branch of $G(\pm i\omega)$ lies above the real axis.   In fact, the relative vertical positions of \(G(i\omega)\) and
\(G(-i\omega)\) determine the stability of the region between the two
branches, since, by Lemma~\ref{localcount} and the strict monotonicity
of \(\Re G(i\omega)\), an upward crossing of \(G(i\omega)\) by
\(1/\gamma\) increases \(N_+(L_\gamma)\) by one, while an upward
crossing of \(G(-i\omega)\) decreases it by one, provided that
\(\Im G(i\omega)\ne0\).

First of all, we will derive the cumulative measure and the corresponding estimates.
    \begin{lemma}\label{cumulative}
    	Assume \(\ell>2\) and \(0<q<\ell-2\). Define the finite
    	signed measure \(\eta\) on \((0,\infty)\) by
    	\begin{align}\label{def-eta}
    		\eta(B):=
    		\int_{\mathbb R\setminus\{s_0,s_1\}}
    		\mathbf1_B(-s)\,d\mu(s),
    		\qquad B\subset(0,\infty)\ \text{Borel}.
    	\end{align}
    	Then
    	\[
    	\operatorname{supp}\eta
    	\subset\{-s_j: 2\leq j\leq J\}\cup[1,\infty),
    	\]
    where
    \[
    -s_j=1-k^2(\ell-2j)^2,
    \qquad
    \eta(\{-s_j\})=\mu_j,
    \qquad j=2,\ldots,J.
    \]
    The set $\{-s_j:2\le j\le J\}$ is empty when $J=1$, and when $J\ge2$, $0<-s_2<\cdots<-s_J<1.$ 
    
    Moreover, the  mass and the first moment satisfy
    \[
    \eta((0,\infty))>0,
    \qquad
    \int_{(-\infty,0)} s\,d\mu(s)>0.
    \]
    In addition, for every $T>0$,
    \begin{equation}\label{cumulativepositive}
    C(T):=\eta((0,T])\ge0,
    \qquad
    \int_{[-T,0)} s\,d\mu(s)\ge0.
    	\end{equation}
    For every \(y\ge0\),
    	\begin{equation}\label{tailkernel}
    		\int_{(0,\infty)}\frac{d\eta(t)}{t^2+y}
    		=\int_0^\infty\frac{2tC(t)}{(t^2+y)^2}\,dt>0.
    	\end{equation}
    \begin{align}\label{tmomenttailestimate}
    \int_{(-\infty,0)}
    \frac{s\,d\mu(s)}{s^2+y}>0,
    \qquad
    \int_{(-\infty,0)}
    \frac{s\,d\mu(s)}{(s^2+y)^2}>0.
    \end{align}
    \end{lemma}
    
    \begin{proof}
    	\textbf{Step 1: Positivity of the total mass and first moment. }  If $2<\ell\le4$, then $J=1$,
so the only discrete spectral weights
are $\mu_0$ and $\mu_1$. Restricting to
    $\mathbb R\setminus\{s_0,s_1\}$ in \eqref{def-eta}
    therefore leaves  the continuous part of $\mu_c$ defined by \eqref{mu-c-B}.  Since $(q-\ell)/2\in(-2,-1)$, one has from  \eqref{signs} that
    $T_q(\xi)>0$ for $\xi>0$. Together with
    $T_{\ell+1}(\xi)>0$ and \eqref{weylmeasure}, this shows
    that $\mu_{\mathrm c}$ is a nonzero positive measure, which implies $\eta((0,\infty))>0$. For $\ell>4$, we prove the same conclusion as follows.
    
    	In light of \eqref{regularsolution}, one has
    	\[
    	\psi_{m_0}=f_\ell,\qquad
    	\psi_{m_1}=(2\ell-1)f_\ell-2(\ell-1)f_{\ell-2},
    	\]
    	where $f_\ell$ is given in \eqref{scalednotation}.
    	Let \(P\) be the orthogonal projection in \(L^2(0,\infty)\)
    	onto \(\operatorname{span}\{f_\ell,f_{\ell-2}\}\).
    	Setting \(X=\sech^2z\), write
    	\[
    	(I-P)f_{\ell+1}=f_{\ell-2}h(X),
    	\qquad
    	X^{3/2}-h(X)\in\operatorname{span}\{1,X\}.
    	\]
    	Since
    	\((I-P)f_{\ell+1}\perp
    	\operatorname{span}\{f_\ell,f_{\ell-2}\}\),
    	we obtain
    	\begin{align}\label{two-orthogonal-condition}
    		\int_0^1
    		h(X)\frac{X^{\ell-3+j}}{\sqrt{1-X}}\,dX=0,
    		\qquad j=0,1.
    	\end{align}
    Moreover, $h''(X)=\frac{3}{4\sqrt X}>0$.
    The orthogonality conditions imply that $h$ must change
    sign at least twice. Indeed, if $h$ did not change sign,
    its weighted integral would be nonzero, contradicting
    \eqref{two-orthogonal-condition}.
    If $h$ changed sign exactly once, at $X_0\in(0,1)$,
    then $(X-X_0)h(X)$ would have a fixed sign and would
    not vanish identically. Its weighted integral would
    be nonzero, whereas
    \eqref{two-orthogonal-condition} gives
    \[
    \int_0^1 (X-X_0)h(X)
    \frac{X^{\ell-3}}{\sqrt{1-X}}\,dX=0,
    \]
    a contradiction. Thus, $h$ has at least two distinct
    zeros in $(0,1)$. Since $h''>0$, Rolle's theorem excludes the  three distinct zeros case. Hence, $h$ has exactly two
    zeros and is negative between them and positive outside.

    Let \(0<a<b<1\) be the two zeros of \(h\).  Let $g$ be the affine interpolation function satisfying
    \[
    g(a)=a^{(q-\ell+2)/2},\qquad
    g(b)=b^{(q-\ell+2)/2},
    \]
    namely,
    \[
    g(X)=\frac{b-X}{b-a}a^{(q-\ell+2)/2}
    +\frac{X-a}{b-a}b^{(q-\ell+2)/2}.
    \]
     Since \((q-\ell+2)/2<0\), strict convexity of $X^{(q-\ell+2)/2}$ yields
    \[
    h(X)\bigl(X^{(q-\ell+2)/2}-g(X)\bigr)>0
    \qquad\text{whenever }h(X)\ne0.
    \]
    	Using orthogonality, we obtain
    	\begin{align}\label{tailtotal}
    		\langle f_q,(I-P)f_{\ell+1}\rangle
    		&=\frac12\int_0^1
    		h(X)\frac{X^{(\ell+q-4)/2}}{\sqrt{1-X}}\,dX
    		\notag\\
    		&=\frac12\int_0^1
    		h(X)\bigl(X^{(q-\ell+2)/2}-g(X)\bigr)
    		\frac{X^{\ell-3}}{\sqrt{1-X}}\,dX>0.
    	\end{align}
    The integrals converge since \((\ell+q-4)/2>-1\)
    	and the singularity at \(X=1\) is integrable.
    	Finally, \(P=E(\{s_0,s_1\})\), where $E$ denotes the spectral measure associated with $L_0$.
    	The definition of \(\eta\) in \eqref{def-eta} therefore yields
    	\[
    	\eta((0,\infty))
    	=c\langle f_q,(I-P)f_{\ell+1}\rangle>0.
    	\]

    We next discuss the positivity of $\int_{(-\infty,0)} s\,d\mu(s)$.     Thanks to Lemma~\ref{moments}, we obtain the total first moment is well-defined and has the finite total variation. Then we set 
    \[
    \begin{aligned}
    \bar h(X)={}&(k^2q^2-1)X^{(q-\ell)/2}
    +k^2(\ell-q)(\ell+q+1)X^{(q-\ell+2)/2}-s_0\frac{\int_0^\infty f_{q+\ell}\,dz}
               {\int_0^\infty f_{2\ell}\,dz}.
    \end{aligned}
    \]
    Invoking the identity  $L_0f_q
    =(k^2q^2-1)f_q
    +k^2(\ell-q)(\ell+q+1)f_{q+2}$ with $f_q$ and $q$ given in \eqref{scalednotation}
    and the self-adjointness of $L_0$, one obtains
    \begin{equation}\label{eq:mixed-first-moment-covariance}
    \int_0^\infty f_{2\ell}\bar h(X)\,dz=0,\quad
    \int_{(-\infty,0)}s\,d\mu(s)
    =c\int_0^\infty f_{2\ell}\bar h(X)\sqrt X\,dz,
    \end{equation}
    where we used  $\mu_0
    =c\frac{\left(\int_0^\infty f_{q+\ell}\,dz\right)
            \left(\int_0^\infty f_{2\ell+1}\,dz\right)}
           {\int_0^\infty f_{2\ell}\,dz}$.
    We claim that $\bar h(1)>0$.  Firstly, evaluating $\bar h$ at $X = 1$ and using \eqref{s0-beta-p-3}, we obtain
    \begin{align}\label{barh1}
    \bar h(1)=s_0\left(
    \frac{3\ell-1-q}{2\ell-1}
    -\frac{\int_0^\infty f_{q+\ell}\,dz}
           {\int_0^\infty f_{2\ell}\,dz}
    \right).
    \end{align}
    By the direct computation, we have the fact that $u\mapsto \sqrt{u-\tfrac12}\int_0^\infty f_u\,dz$
    is strictly increasing for $u>1/2$. Thus
    \begin{equation}\label{eq:mixed-strict-Wallis}
    \frac{\int_0^\infty f_{q+\ell}\,dz}
         {\int_0^\infty f_{2\ell}\,dz}
    <
    \sqrt{\frac{2\ell-1/2}{q+\ell-1/2}}
    <
    \frac{3\ell-1-q}{2\ell-1}.
    \end{equation}
     It then follows that our claim holds thanks to \eqref{barh1}.   
Moreover,
\[
\begin{aligned}
\bar h'(X)=X^{(q-\ell-2)/2}\biggl[
(k^2q^2-1)\frac{q-\ell}{2}+k^2(\ell-q)(\ell+q+1)
\frac{q-\ell+2}{2}X
\biggr].
\end{aligned}
\]
The expression in brackets is strictly decreasing
and positive at $X=0$. Thus $\bar h$ first
increases and then possibly decreases.
Since $\lim_{X\downarrow0}\bar h(X)=-\infty$
and $\bar h(1)>0$, it has exactly one zero
$\bar a\in(0,1)$, with
\[
\bar h(X)<0\quad(0<X<\bar a),\qquad
\bar h(X)>0\quad(\bar a<X<1).
\]

    Therefore, we get
    \[
    \int_0^\infty f_{2\ell}{\bar h}(X)\sqrt X\,dz
    =
    \int_0^\infty f_{2\ell}\bar h(X)(\sqrt X-\sqrt {\bar a})\,dz
    >0.
    \]
     This proves positivity of $\int_{(-\infty,0)}s\,d\mu(s)$ in
    \eqref{eq:mixed-first-moment-covariance}.

   \textbf{Step 2:  Positivity of the cumulative mass and the first moment.}
    By \eqref{atommass}, \eqref{boundtransform-explicit},
    and \eqref{Ni-normalization}, for $0\le j<J$ with
    $\mu_j\ne0$, we have
    \begin{equation}\label{atomratio}
    	\frac{\mu_{j+1}}{\mu_j}
    	=\frac{m_j-2}{m_j}
    	\frac{(j+1/2)^2}{j+1}
    	\frac{\ell-j}{(\ell-j-1/2)^2}
    	\frac{q-m_j}{q+m_j-2}.
    \end{equation}
    Since $m_j>2$ and $q>0$, one gets
    \[
    \begin{aligned}
    	0\le
    	\frac{m_j-2}{m_j}\frac{|q-m_j|}{q+m_j-2}
    	&<
    	\frac{m_j-2}{m_j}\frac{q+m_j}{q+m_j-2}=1-\frac{2q}{m_j(q+m_j-2)}<1.
    \end{aligned}
    \]
   Also, $\ell-j>j+2$, and
$u/(u-1/2)^2$ is strictly decreasing
for $u>1/2$. Hence,
    \begin{align}\label{remaining-factor-second-factor}
    	0<
    	\frac{(j+1/2)^2}{j+1}
    	\frac{\ell-j}{(\ell-j-1/2)^2}
    	<
    	\frac{(j+1/2)^2(j+2)}
    	{(j+1)(j+3/2)^2}
    	<1,
    \end{align}    
Combining these estimates with \eqref{atomratio},
we obtain
\begin{equation}\label{consecute-weight}
|\mu_{j+1}|<|\mu_j|,
\qquad 0\le j<J,\quad \mu_j\ne0.
\end{equation}

   For $1\le j<J$ with $\mu_j\mu_{j+1}<0$,
we have $0<q<m_j$. Using \eqref{atomratio} we obtain
    \begin{align}\label{collect1-ratio}
    \left|\frac{s_{j+1}\mu_{j+1}}{s_j\mu_j}\right|
    &=
    \frac{s_{j+1}}{s_j}
    \frac{m_j-2}{m_j}
    \frac{m_j-q}{m_j+q-2}
    \frac{(j+1/2)^2}{j+1}
    \frac{\ell-j}{(\ell-j-1/2)^2}.\\
    &<
	\frac{(2j+1)^4(j+2)}
	{(j+1)(2j-1)(2j+3)^3}<1
    \end{align}
Thus,
\(|s_{j+1}\mu_{j+1}|<|s_j\mu_j|\)
whenever \(\mu_j\mu_{j+1}<0.\) for the alternating pairs.
     
    Set $n=\lceil(\ell-q)/2\rceil$. Since
    $1<(\ell-q)/2<\ell/2$, we have $2\le n\le J+1$.
    By \eqref{signs},
    \begin{equation}\label{pattern}
    \begin{aligned}
    \operatorname{sgn}\mu_j
    &=(-1)^j, & 2\le j\le\min(n,J),\\
    \operatorname{sgn}(s_j\mu_j)
    &=(-1)^{j+1}, & 1\le j\le\min(n,J).
    \end{aligned}
    \end{equation}
    
If $(\ell-q)/2\notin\mathbb Z$, then by
\eqref{signs}, the weights $\mu_j$ with
$n<j\le J$ and the continuous part $\mu_{\mathrm c}$
have sign $(-1)^n$, while $s_j\mu_j$ and
$s\,d\mu_{\mathrm c}(s)$ have sign $(-1)^{n+1}$.
If $(\ell-q)/2=n\in\mathbb Z$, then
$\mu_j=0$ for $j>n$ and $\mu_{\mathrm c}=0$.

 For $J\ge2$, the weights
$\mu_2,\ldots,\mu_{\min(n,J)}$ alternate in sign,
starting with $\mu_2>0$. By \eqref{consecute-weight},
pairing terms gives
\[
\sum_{i=2}^j\mu_i>0,
\qquad 2\le j\le\min(n,J).
\]
Thus, $C(T)>0$ for
$-s_2\le T\le-s_{\min(n,J)}$.
By the sign properties,
$C$ is monotone on
$[-s_{\min(n,J)},\infty)$.
Since $\eta((0,\infty))>0$ by Step 1,
we conclude that $C(T)>0$ for $T\ge-s_2$.
Also, $C(T)=0$ for $0<T<-s_2$.
When $J=1$, $\eta$ is a positive measure
by Step 1, so $C(T)\ge0$ for every $T>0$.

Similarly, the weights
$s_1\mu_1,\ldots,s_{\min(n,J)}\mu_{\min(n,J)}$
alternate in sign, starting with $s_1\mu_1>0$.
By the weighted ratio estimate above,
\[
\sum_{i=1}^j s_i\mu_i>0,
\qquad 1\le j\le\min(n,J).
\]
The remaining first moment contributions have
a fixed sign or vanish. Hence,
$\int_{[-T,0)}s\,d\mu(s)$ is monotone for
$T\ge-s_{\min(n,J)}$.
Together with the positivity of the total
first moment from Step 1, this yields
\[
\int_{[-T,0)}s\,d\mu(s)>0,
\qquad T\ge-s_1.
\]
The integral vanishes for $0<T<-s_1$.
This proves \eqref{cumulativepositive}.
    
    \textbf{Step 3: Positivity of the integral.}
    Fix \(y\ge0\). Since $C(t)=0$ for all sufficiently small $t>0$ and \(|C(t)|\le\|\eta\|_{\mathrm{TV}}\), we have
    \[
    \lim_{t\downarrow0}\frac{C(t)}{t^2+y}=0,
    \qquad
    \lim_{t\to\infty}\frac{C(t)}{t^2+y}=0.
    \]
    Then,  by using \(C(t)=\eta((0,t])\), we obtain from the integration by parts  that 
    \[
    \begin{aligned}
    	\int_{(0,\infty)}\frac{d\eta(t)}{t^2+y}
    	&=
    	\left[\frac{C(t)}{t^2+y}\right]_{0^+}^{\infty}
    	-\int_0^\infty C(t)
    	\frac{d}{dt}\left(\frac1{t^2+y}\right)\,dt\\
    	&=
    	\int_0^\infty\frac{2tC(t)}{(t^2+y)^2}\,dt>0,
    \end{aligned}
    \]
    where for $J\ge2$, the strict inequality follows from Step~2,
    which gives $C(t)>0$ for $t\ge -s_2$. For $J=1$,
    $\eta$ is a nonzero positive measure by Step~1, so $C(t)$
    is positive for all sufficiently large $t$.
    This proves \eqref{tailkernel}.

    For \(j=1,2\), we similarly apply the integration by parts to obtain
    \[
    \int_{(-\infty,0)} \frac{s\,d\mu(s)}{(s^2+y)^j}
    =
    \int_0^\infty \frac{2jt}{(t^2+y)^{j+1}}
    \left(\int_{[-t,0)} s\,d\mu(s)\right)\,dt.
    \]
    The strict inequality follows from Step~2, which gives $\int_{[-t,0)} s\,d\mu(s)>0$ for  $t\ge -s_1.$
    It completes the proof of \eqref{tmomenttailestimate}.
    
    \end{proof}
    
    In the following corollary, we show that for $0<q\le\ell$,
    $\Re G(i\omega)$ decreases strictly from $1$ to $0$ as
    $\omega$ increases from $0$ to $\infty$. The proof of this corollary is based on
    Lemma~\ref{cumulative} in the mixed-sign
    case. Later on, we use this
    monotonicity to parametrize each curve $G(\pm i\omega)$
    by its real part in the proof of Theorem~\ref{main}.
    \begin{corollary}\label{complex-realpart}
    If $0<q\le\ell$, equivalently $1<r\le(p+3)/2$, then
    \[
    \Re G(i\omega)>0\quad(\omega\ge0),
    \qquad
    \frac{d}{d\omega}\Re G(i\omega)<0\quad(\omega>0).
    \]
    Moreover,
    \[
    G(0)=1,\qquad
    \lim_{\omega\to\infty}\Re G(i\omega)=0.
    \]
    In particular, $\omega\mapsto G(i\omega)$ is injective on
    $[0,\infty)$ and $G'(i\omega)\ne0$ for $\omega>0$.
    \end{corollary}
    \begin{proof}
    For $\omega\ge0$, set $y=\omega^2$.  By using the spectral representation given in \eqref{tv-bound-G-spectral}, we have
    
    \begin{align}\label{Re-G-derivative}
    \begin{aligned}
    \Re G(i\sqrt y)
    &=\frac{s_0\mu_0}{s_0^2+y}
    +\int_{(-\infty,0)}\frac{s\,d\mu(s)}{s^2+y},\\
    -\frac{d}{dy}\Re G(i\sqrt y)
    &=\frac{s_0\mu_0}{(s_0^2+y)^2}
    +\int_{(-\infty,0)}\frac{s\,d\mu(s)}{(s^2+y)^2}.
    \end{aligned}
    \end{align}
    If $0<q<\ell-2$ (so $\ell>2$), We invoke Lemma~\ref{cumulative} to obtain
    \begin{align}\label{both-nonnegative-holds}
    \int_{(-\infty,0)}\frac{s\,d\mu(s)}{(s^2+y)^j}>0,
    \qquad j=1,2,\quad y\ge0.
    \end{align}
    If $0<q<\ell$ and $q\ge\ell-2$, then
    $\mu|_{(-\infty,0)}\le0$ thanks to \eqref{signs}. Moreover, since $s<0$ on $(-\infty,0)$, $ \int_{(-\infty,0)}\frac{s\,d\mu(s)}{(s^2+y)^j}\geq 0$ for 
    $j=1,2$ and $y\ge0.$
    If $q=\ell$,  one finds the corresponding two integrals  vanish.
    In all three cases, the two integrals over $(-\infty,0)$
    are nonnegative. Since $s_0\mu_0>0$, both
    $\Re G(i\sqrt y)$ and $-\frac{d}{dy}\Re G(i\sqrt y)$
    are strictly positive for $y\ge0$.
    
    In addition, for $\omega>0$, we use 
    the chain rule to get
    \[
    -\Im G'(i\omega)
    =\frac{d}{d\omega}\Re G(i\omega)
    =2\omega\left.\frac{d}{dy}\Re G(i\sqrt y)\right|_{y=\omega^2}
    <0,
    \]
    which proves the  injectivity of $G$ on $[0,\infty)$ and
    $G'(i\omega)\ne0$ for $\omega>0$. In light of Lemma~\ref{moments}, one has
    $G(0)=1$. As $\omega\to\infty$, we have $y\to\infty$, 
    so $\Re G(i\omega)\to0$ by \eqref{Re-G-derivative} and the  dominated convergence theorem, since
    $s\,d\mu(s)$ has finite total variation.
    \end{proof}

Next, we analyze the sign of $H(y)$ defined in
\eqref{def-H-y-s2-y}, which determines the sign of
$\Im G(i\omega)$ through the identity
$\Im G(i\omega)=\omega H(\omega^2)$ for $\omega>0$.

      \begin{lemma}\label{frequency}
{
For every \(\omega\in\mathbb R\), \(G(i\omega)\ne0\).
For every \(\gamma\in\mathbb C\) with \(\Re\gamma\ge0\),
\begin{equation}\label{complex-axis-equivalence}
\sigma(L_\gamma)\cap i\mathbb R
=
\left\{i\omega:\omega\in\mathbb R,\ 
\gamma=\frac{1}{G(i\omega)}\right\}.
\end{equation}
}
The following statements hold.
\begin{enumerate}
\item If \(r\ge(p-1)/2\), then \(H(y)>0\) for every \(y>0\), so \(G(i\omega)\) lies above \(G(-i\omega)\).
The curve
\(\{1/G(i\omega):\omega\in\mathbb R\}\) intersects the real
axis only at \(\gamma=1\), corresponding to \(\omega=0\).

\item If \(1<r<(p-1)/2\), there exists a unique \(\omega_H>0\)
such that \(H(\omega_H^2)=0\). Moreover,
\[
0<\Re G(i\omega_H)<1,
\qquad
\Re G'(i\omega_H)>0.
\]
Define \(\gamma_H:=1/\Re G(i\omega_H)>1\), then $\Im G(i\omega)<0$ for $0<\omega<\omega_H$ and $\Im G(i\omega)>0$ for $\omega>\omega_H$.  In addition, \(G(-i\omega)\) lies above \(G(i\omega)\) for
\(0<\omega<\omega_H\), and \(G(i\omega)\) lies above
\(G(-i\omega)\) for \(\omega>\omega_H\).
At \(\omega=\omega_H\), both curves pass through \(1/\gamma_H\).

{Also, the curve \(\{1/G(i\omega):\omega\in\mathbb R\}\)
intersects the real axis at \(1\) and \(\gamma_H\), with
\(\gamma_H\) attained precisely at \(\omega=\pm\omega_H\).}
The eigenvalues \(\pm i\omega_H\) of \(L_{\gamma_H}\) are
algebraically simple, and their local analytic branches
\(\lambda_\pm(\gamma)\)
{for complex \(\gamma\) near \(\gamma_H\)}
satisfy
\[
\lambda_\pm(\gamma_H)=\pm i\omega_H,
\qquad
\Re\lambda_\pm'(\gamma_H)
=
-\frac{\Re G'(i\omega_H)}
{\gamma_H^2|G'(i\omega_H)|^2}<0.
\]
\end{enumerate}
\end{lemma}

        \begin{proof}
    		We use \(H(0)=G'(0)\) and consider the following three cases.

    \textbf{Case 1: $q\ge\ell$.}
    By \eqref{signs}, $\mu$ is nonnegative and $\mu_0>0$,
    so $H(y)>0$ for all $y\ge0$.

    		\textbf{Case 2:
    			\(\ell-2\le q<\ell\).}
    		Since \((q-\ell)/2\in[-1,0)\), \eqref{signs} gives that
    		\(\mu|_{(-\infty,0)}\le0\).
    		If \(\ell>2\), then \(\mu_1<0\);
    		if \(1<\ell\le2\), then \(T_q(\xi)<0\) for all
    		\(\xi>0\). Thus, \(\mu|_{(-\infty,0)}\ne0\).
    		
    		Since \(0\in\rho(L_0)\), the support of \(\mu\) is
    		bounded away from zero.  In light of \eqref{spectral-representation-G}, we obtain
    		\begin{equation}\label{stieltjes}
    			G(\lambda)
    			=\frac{\mu_0}{s_0-\lambda}
    			+\int_{(-\infty,0)}\frac{d\mu(s)}{s-\lambda},
    			\qquad \lambda\in\rho(L_0),
    		\end{equation}
    		and
    		\[
    		H(y)
    		=\frac{\mu_0}{s_0^2+y}
    		+\int_{(-\infty,0)}\frac{d\mu(s)}{s^2+y},
    		\qquad y\ge0,
    		\]  
            where $H$ is defined by \eqref{def-H-y-s2-y}.
    	By \eqref{stieltjes} and Lemma~\ref{moments},
    	\begin{equation}\label{hrelations}
    		H(0)=G'(0)=\frac1{p-1}-\frac1{2r},
    		\qquad
    		\lim_{y\to\infty}(s_0^2+y)H(y)
    		=\mu(\mathbb R)>0.
    	\end{equation}
    Then, we differentiate to get
    		\begin{equation}\label{hshape}
    			\begin{aligned}
    				\frac{d}{dy}\bigl((s_0^2+y)H(y)\bigr)
    				&=\int_{(-\infty,0)}
    				\frac{s^2-s_0^2}{(s^2+y)^2}\,d\mu(s),\\
    				\frac{d}{dy}\left[
    				(s_0^2+y)^2
    				\frac{d}{dy}\bigl((s_0^2+y)H(y)\bigr)
    				\right]
    				&=2(s_0^2+y)
    				\int_{(-\infty,0)}
    				\frac{(s^2-s_0^2)^2}{(s^2+y)^3}\,d\mu(s)<0,
    			\end{aligned}
    		\end{equation}
    		where the last inequality is strict since
    		\(\mu|_{(-\infty,0)}\) is not supported only at
    		\(s=-s_0\).   By using \eqref{signs} , we find if \(\ell>2\), then
    		\[
    		\mu_1<0,\qquad -s_0<s_1<0	. 
    		\]
    	 If \(1<\ell\le2\), the continuous part of \(\mu\)
    		is nonzero.

    By \eqref{hshape}, the derivative of
    $(s_0^2+y)H(y)$ has at most one zero in $(0,\infty)$.
    If it vanishes, its sign changes from positive
    to negative. If $r\ge(p-1)/2$, then $H(0)=G'(0)\ge0$.
    Together with the positive limit in
    \eqref{hrelations}, this implies $H(y)>0$
    for every $y>0$.
    
    If $r<(p-1)/2$, then $H(0)<0$, while
    the limit in \eqref{hrelations} is positive.
    Thus, $H$ has a positive zero. At any such zero,
    the derivative of $(s_0^2+y)H(y)$ must be positive;
    otherwise, \eqref{hshape} would imply a decreasing
    function, contradicting its positive limit.
    Since this derivative can change sign at most once,
    the positive zero is unique. Denoting it by
    $\omega_H^2$, we obtain
    
    \[
    (s_0^2+\omega_H^2)H'(\omega_H^2)
    =
    \left.
    \frac{d}{dy}\bigl((s_0^2+y)H(y)\bigr)
    \right|_{y=\omega_H^2}>0.
    \]
    By \eqref{stieltjes}, for every $\omega>0$,
    \[
    \begin{aligned}
    	0<\Re G(i\omega)
    	&=\frac{\mu_0 s_0}{s_0^2+\omega^2}
    	+\int_{(-\infty,0)}
    	\frac{s}{s^2+\omega^2}\,d\mu(s)\\
    	&<\frac{\mu_0}{s_0}
    	+\int_{(-\infty,0)}\frac{d\mu(s)}{s}
    	=G(0)=1,
    \end{aligned}
    \]
   where we used $\mu|_{(-\infty,0)}\le0$
by \eqref{signs}.  Finally, differentiating
    		\(\Im G(i\omega)=\omega H(\omega^2)\) and using
    		\(H(\omega_H^2)=0\), we find
    		\begin{equation}\label{crossingderivative}
    			\Re G'(i\omega_H)
    			=
    			\left.\frac{d}{d\omega}\Im G(i\omega)\right|_{\omega=\omega_H}
    			=
    			2\omega_H^2H'(\omega_H^2)>0.
    		\end{equation}
    		\textbf{Case 3: \(0<q<\ell-2\).}
    		Here \(\ell>2\), so \(1<p<3\). Consequently, $r>1>\frac{p-1}{2}$, and  $G'(0)>0.$
    		Moreover, \(\mu_1<0\)  and a direct calculation gives
    		\begin{equation}\label{firstlocations}
    			-s_1=\frac{2\ell-3}{(\ell-1)^2},
    			\qquad
    			s_0=\frac{2\ell-1}{(\ell-1)^2},
    			\qquad
    			0<-s_1<\min\{1,s_0\}.
    		\end{equation}
    		\textit{If \(2<\ell\le4\),}
    		the only negative discrete eigenvalue of \(L_0\) on the even
    		subspace is \(s_1\). Since
    		\((q-\ell)/2\in(-2,-1)\), \eqref{signs} shows that
    		the continuous part \(\mu_{\mathrm c}\)  is a positive measure.
    		Thus,
    		\[
    		H(y)
    		=\frac{\mu_0}{s_0^2+y}
    		+\frac{\mu_1}{s_1^2+y}
    		+\int_{(-\infty,-1]}
    		\frac{d\mu_{\mathrm c}(s)}{s^2+y}.
    		\]
    		It then follows that
    		\[
    		\begin{aligned}
    			\frac{d}{dy}\bigl((s_1^2+y)H(y)\bigr)
    			&=\frac{\mu_0(s_0^2-s_1^2)}{(s_0^2+y)^2}\\
    			&\quad+\int_{(-\infty,-1]}
    			\frac{s^2-s_1^2}{(s^2+y)^2}\,d\mu_{\mathrm c}(s)>0,
    		\end{aligned}
    		\]
    		where we used \(\mu_0>0\), \(s_0>-s_1\)  and
    		\(s^2\ge1>s_1^2\) on the continuous spectrum.
    		Since $\left.(s_1^2+y)H(y)\right|_{y=0}
    		=s_1^2G'(0)>0,$
    		we therefore conclude that \(H(y)>0\) for every \(y\ge0\).
    		
    		\textit{If \(\ell>4\),}
    		we invoke Lemma~\ref{cumulative} to get
    		\[
    		H(y)
    		=\frac{\mu_0}{s_0^2+y}
    		+\frac{\mu_1}{s_1^2+y}
    		+\int_{(0,\infty)}\frac{d\eta(t)}{t^2+y},
    		\]
    		where for $y\ge0$, $\int_{(0,\infty)}\frac{d\eta(t)}{t^2+y}>0.$  By \eqref{atomratio} with $j=0$ and
    $s_0>-s_1$, for every $y\ge0$,
    \begin{equation}\label{pairbound}
    	\begin{aligned}
    		0<
    		-\frac{\mu_1(s_0^2+y)}{\mu_0(s_1^2+y)}
    		&\le-\frac{\mu_1s_0^2}{\mu_0s_1^2}=\frac{(\ell-2)(\ell-q)}
    		{(2\ell-3)^2(q+\ell-2)}
    		<\frac{\ell}{(2\ell-3)^2}<1,
    	\end{aligned}
    \end{equation}
    where we used $q>0$ and $\ell>4$.
    Together with Lemma~\ref{cumulative}, this yields
    \begin{equation}\label{mixedbound}
    	H(y)>
    	\frac{\mu_0}{s_0^2+y}
    	+\frac{\mu_1}{s_1^2+y}>0.
    \end{equation}

   Since $\Im G(i\omega)=\omega H(\omega^2)$ for $\omega>0$,
  $\omega_H$ is the unique positive
frequency for which $G(i\omega)$ is real  in the case $1<r<(p-1)/2$.
At this frequency, $G(i\omega_H)=\Re G(i\omega_H)>0$.
Together with $G(0)=1$ shown in Lemma \ref{moments} and
$G(-i\omega)=\overline{G(i\omega)}$, this shows that
$G(i\omega)\ne0$ for all $\omega\in\mathbb R$. Thus, the curve $\{1/G(i\omega):\omega\in\mathbb R\}$
intersects the real axis only at $1$ if $r\ge(p-1)/2$,
and at $1$ and $\gamma_H$ if $1<r<(p-1)/2$.

Now, we determine which branch of $G(\pm i\omega)$ lies above the real axis. For \(r\ge(p-1)/2\), the positivity of \(H\) yields
\(\Im G(i\omega)>0\) for all \(\omega>0\). For \(1<r<(p-1)/2\), by using 
\eqref{crossingderivative}, one has that \(\Im G(i\omega)\) changes
from negative to positive at \(\omega_H\). By continuity and
uniqueness of the positive zero, we obtain
\[
\Im G(i\omega)<0 \quad (0<\omega<\omega_H),
\qquad
\Im G(i\omega)>0 \quad (\omega>\omega_H).
\]
Then, the ordering of \(G(i\omega)\) and \(G(-i\omega)\) stated in the lemma follows from
\(\Im G(-i\omega)=-\Im G(i\omega)\).

By Lemma~\ref{local}, $i\mathbb R\subset\rho(L_0)$.
Hence, Lemma~\ref{rankone} gives us
\[
i\omega\in\sigma(L_\gamma)
\quad\Longleftrightarrow\quad
1-\gamma G(i\omega)=0
\quad\Longleftrightarrow\quad
\gamma=\frac{1}{G(i\omega)},
\qquad \omega\in\mathbb R.
\]
This proves \eqref{complex-axis-equivalence}.
    
  If $1 <r<(p-1)/2$, then  thanks to \eqref{crossingderivative},
    \(\Re G'(i\omega_H)>0\), and hence
    \[
    d_{\gamma_H}'(\pm i\omega_H)
    =-\gamma_HG'(\pm i\omega_H)\ne0.
    \]
    Thus, \(\pm i\omega_H\) are simple zeros of
    \(d_{\gamma_H}\), and Lemma~\ref{rankone} shows that
    both eigenvalues are algebraically simple.  The implicit function theorem further yields local analytic
    branches $\lambda_\pm(\gamma)$ satisfying
    \[
    \lambda_\pm(\gamma_H)=\pm i\omega_H,
    \qquad
    1-\gamma G(\lambda_\pm(\gamma))=0.
    \]
    		Differentiating at \(\gamma=\gamma_H\), one gets $\lambda_\pm'(\gamma_H)
    		=-\frac{1}{\gamma_H^2G'(\pm i\omega_H)}.$
    		Therefore,
    		\begin{equation}\label{direction}
    			\Re\lambda_\pm'(\gamma_H)
    			=-\frac{\Re G'(i\omega_H)}
    			{\gamma_H^2|G'(i\omega_H)|^2}<0.
    		\end{equation}
    	Thus, both eigenvalues cross the imaginary axis transversely from right
to left along any differentiable path through \(\gamma_H\) whose
tangent direction is sufficiently close to the positive real direction.
    	\end{proof}
    
    \section{Spectrum analysis of NLEP}
In this section, we prove Theorem~\ref{main} with the aid of 
Lemmas~\ref{cumulative} and~\ref{frequency} and
Corollary~\ref{complex-realpart}. The curve
$\gamma=1/G(i\omega)$, $\omega\in\mathbb R$, consists of $\gamma$ for which $L_\gamma$ has an eigenvalue on the imaginary
axis. Hence, we use the curves $G(\pm i\omega)$ to count the eigenvalues in the open
right half-plane, both on the curve and in the regions it separates. The regimes is determinedn by the sign of  $\mu$ and
the behavior of $H$. When $r>(p+3)/2$,
$\mu$ is nonnegative. For $r<(p+3)/2$
it has mixed signs. In the latter case, $\Re G(i\omega)$ is
strictly decreasing. The sign of
$H(0)=G'(0)$ determines whether $H$
has a positive zero, with the threshold
at $r=(p-1)/2$.

To be more precise, at $r=(p+3)/2$, as $\omega$ increases from $0$, $G(i\omega)$ traces
the upper semicircle from $1$ toward $0$. Taking reciprocals maps
this branch and its conjugate to the vertical line $\Re\gamma=1$.
By Lemma~\ref{moments}, $G(0)=1$, so the conjugate branches
$G(\pm i\omega)$ meet at $1$ when $\omega=0$. This point
corresponds to $\gamma=1$, where $L_1$ has a zero eigenvalue.
When $r<(p+3)/2$, Lemma~\ref{frequency} shows that the branches
have a second real intersection precisely when $1<r<(p-1)/2$.

Moreover, at the threshold $r=(p-1)/2>1$, we have $G'(0)=0$.
The following lemma shows that the zero eigenvalue of $L_1$
has algebraic multiplicity two and describes how it splits as
$\gamma\to1$ in $\mathbb C$. We use the expansions
to determine $N_+(L_1)$ at this threshold.
    	\begin{lemma}\label{critical}
    		Assume \(r=(p-1)/2>1\). Then \(p>3\),
    		\(G\) satisfies \eqref{stieltjes}, and
    		\[
    		G''(0)>0,\qquad G'''(0)<0.
    		\]
    		Moreover,
    		\[
    		d_1(0)=d_1'(0)=0,
    		\qquad
    		d_1''(0)=-G''(0)<0.
    		\]
    		Thus, zero is an eigenvalue of \(L_1\) of algebraic
    		multiplicity two.  As $\gamma\rightarrow 1$ in $\C$, the two eigenvalues $\lambda_\pm(\gamma)$ of $L_{\gamma}$ near zero satisfy \eqref{criticalexpansion}.
    	\end{lemma}
    	\begin{proof}
    		Since \(r=(p-1)/2>1\), we have
    		\[
    		p>3,\qquad
    		1<\ell<2,\qquad
    		q=\frac{p-3}{p-1}>0,\qquad
    		\frac{q-\ell}{2}=-\frac{2}{p-1}\in(-1,0).
    		\]
    	By \eqref{signs}, $\mu|_{(-\infty,0)}$ is
    	nonpositive and has a nonzero continuous part.
    	Hence, the representation \eqref{stieltjes} applies.
    		In light of \(G'(0)=0\), we invoke Lemma~\ref{moments} to find
    		\[
    		0<-\int_{(-\infty,0)}d\mu(s)<\mu_0,
    		\qquad
    		-\int_{(-\infty,0)}\frac{d\mu(s)}{s^2}
    		=\frac{\mu_0}{s_0^2}.
    		\]
    		Applying Cauchy--Schwarz inequality to the positive measure
    		\(-\mu|_{(-\infty,0)}\), we further obtain
    		\begin{equation}\label{criticalCS}
    			\begin{aligned}
    				-\int_{(-\infty,0)}\frac{d\mu(s)}{s^4}
    				&\ge
    				\frac{\left(
    					-\int_{(-\infty,0)}s^{-2}\,d\mu(s)
    					\right)^2}
    				{-\int_{(-\infty,0)}d\mu(s)}
    				=
    				\frac{\mu_0^2}
    				{s_0^4\left(-\int_{(-\infty,0)}d\mu(s)\right)}
    				>\frac{\mu_0}{s_0^4}.
    			\end{aligned}
    		\end{equation}
    	Differentiating \eqref{stieltjes}, we get
    		\begin{equation}\label{criticalmoments}
    			\begin{aligned}
    				G''(0)
    				&=2\frac{\mu_0}{s_0^3}
    				+2\int_{(-\infty,0)}\frac{d\mu(s)}{s^3}>0,\\
    				G'''(0)
    				&=6\left(
    				\frac{\mu_0}{s_0^4}
    				+\int_{(-\infty,0)}\frac{d\mu(s)}{s^4}
    				\right)<0,
    			\end{aligned}
    		\end{equation}
    		where the first inequality holds by using  \(s^3<0\) on \((-\infty,0)\);
    		the second follows from \eqref{criticalCS}.

    		Since \(G(0)=1\) and \(G'(0)=0\) as shown in Lemma \ref{moments},
    		\[
    		d_1(\lambda)
    		=-\frac{G''(0)}2\lambda^2
    		-\frac{G'''(0)}6\lambda^3+O(\lambda^4).
    		\]
    	Moreover, by using \(G''(0)>0\) shown in \eqref{criticalmoments}, the function \(d_1\) has a zero of order exactly two at \(\lambda=0\). Lemma~\ref{rankone} therefore
    		shows that 0 is an eigenvalue of \(L_1\) of algebraic
    		multiplicity two. By Lemma~\ref{moments}, $\ker L_1=\operatorname{span}\{w\}.$
    Consequently, Lemma~\ref{localcount} implies that, for a
    	sufficiently small disk $D$ centered at zero, $n_D(L_\gamma)=2$,
    	whenever $\gamma$ is sufficiently close to $1$. 
    		
    	Set $\gamma=1+t^2$ with $t\in\mathbb C$ near zero  and
    	$\lambda=tv$. The equation $d_\gamma(\lambda)=1-\gamma G(\lambda)=0$
    	becomes
    		\begin{align}\label{take-deriative-ratio}
    			\frac{1-(1+t^2)G(tv)}{t^2}
    			=-1-\frac{G''(0)}2v^2
    			-\frac{G'''(0)}6tv^3+O(t^2)=0.
    		\end{align}
            Since \(G(0)=1\) and \(G'(0)=0\), the quotient extends
    holomorphically to \(t=0\). At \(t=0\), its two roots are $v=\pm i\sqrt{\frac2{G''(0)}}.$
    Both roots are simple, since the derivative with
    respect to $v$ is $-G''(0)v\ne0$.  Invoking the implicit function theorem, one obtains a unique analytic solution $v(t)$ near each root.  For sufficiently small complex $t\not=0$, the two solutions
    	$v(t)$ give two distinct eigenvalues
    	$\lambda=tv(t)$ of $L_\gamma$ in $D$.
    	By Lemma~\ref{localcount}, these are the only
    	eigenvalues in $D$.
    	 In addition, we differentiate in \eqref{take-deriative-ratio} to arrive at
    		\[
    		-G''(0)v(0)v'(0)
    		-\frac{G'''(0)}6v(0)^3=0.
    		\]
    		Therefore,
    		\[
    		v(0)=\pm i\sqrt{\frac2{G''(0)}},
    		\qquad
    		v'(0)=\frac{G'''(0)}{3(G''(0))^2}.
    		\]
    		Using \(\lambda=tv(t)\) and \(\gamma=1+t^2\), we prove
    \eqref{criticalexpansion}.
    	\end{proof}
    	
    	Now we are ready to state the proof of our main result.
    	\begin{proof}[Proof of Theorem \ref{main}]
    By Lemma~\ref{local}, $L_0$ has exactly one
    positive eigenvalue and no spectrum on the imaginary axis .
    By Lemmas~\ref{moments} and~\ref{frequency},
    the spectrum of $L_\gamma$ on the imaginary axis is empty
    for $\gamma>0$, except at $\gamma=1$ and, when
    $1<r<(p-1)/2$, at $\gamma=\gamma_H$.
    Hence, Lemma~\ref{rankone} gives
    		\[
    		N_+(L_\gamma)=N_+(L_0)=1,
    		\qquad 0\le\gamma<1.
    		\]	
    At $\gamma=1$, we have
    		\[
    		\sigma(L_1)\cap i\mathbb R=\{0\},
    		\qquad
    		\ker L_1=\operatorname{span}\{w\}.
    		\] Choose a disk \(D\) centered at zero as in
    		Lemma~\ref{localcount}, with \(\gamma_*=1\).
    		For \(\gamma\) sufficiently close to \(1\),
    		\[
    		n_D(L_\gamma)=
    		\begin{cases}
    			1,&r\ne(p-1)/2,\\
    			2,&r=(p-1)/2.
    		\end{cases}
    		\]
    		If \(r\ne(p-1)/2\), the unique eigenvalue  of \(L_\gamma\) in \(D\)
    		satisfies \eqref{simplebranch}.
    		If \(r=(p-1)/2\), the two eigenvalues of \(L_\gamma\) in \(D\)
    		are described in Lemma~\ref{critical}.
    		
    		By Lemma~\ref{localcount},
    		\[
    		N_+(L_\gamma)
    		=
    		N_+(L_1)
    		+n_{D\cap\{\lambda\in\mathbb C:\Re\lambda>0\}}(L_\gamma)
    		\]
    		for \(\gamma\) sufficiently close to \(1\).
    
    \noindent\textbf{The case $r <(p+3)/2$.}  By Corollary~\ref{complex-realpart},	$\Re G(i\omega)$
decreases strictly from $1$ to $0$ as $\omega$ increases.  {
By \eqref{spectral-representation-G}, $G(\pm i\omega)\to0$ as
$\omega\to\infty$, so both graphs $G(\pm i\omega)$ are bounded. Since $G(\pm i\omega)$ are bounded graphs over
$(0,1)$, the exterior of their enclosed regions
is path connected in the right half-plane 
with the origin excluded. Indeed, every exterior point can be joined to
$\Re(1/\gamma)>1$ by moving vertically away
from both curves and then horizontally,
without crossing either curve.

Since $N_+(L_\gamma)=1$ for $0<\gamma<1$,
which corresponds to $1/\gamma>1$,
the path connectedness of the exterior,
\eqref{complex-axis-equivalence},
and Lemma~\ref{rankone} imply that
$N_+(L_\gamma)=1$ throughout the exterior.

}
At $1/\gamma=G(i\omega)$ with $\omega>0$
and $\Im G(i\omega)\ne0$,
Corollary~\ref{complex-realpart} ensures that
$i\omega$ is the only  eigenvalue on imaginary axis.
Since $G'(i\omega)\ne0$, it is algebraically
simple by Lemma~\ref{rankone}.
Fixing $\Re(1/\gamma)$ and differentiating
$G(\lambda)=1/\gamma$ with respect to
$\Im(1/\gamma)$, we obtain
\begin{align}\label{direction-crossing-eigenvalue}
{\frac{d}{d(\Im (1/\gamma))}\Re\lambda}
=
\Re\!\left(\frac{i}{G'(i\omega)}\right)
=
-\frac{\dfrac{d}{d\omega}\Re G(i\omega)}
       {|G'(i\omega)|^2}
>0.
\end{align}
For sufficiently small  $t\in \mathbb{R}$, let \(\frac1{\gamma(t)}=G(i\omega)+it.\)
The exterior corresponds to $t>0$
if $\Im G(i\omega)>0$, and to $t<0$
if $\Im G(i\omega)<0$.
By \eqref{direction-crossing-eigenvalue},
and Lemma~\ref{localcount}, we obtain

\begin{align}\label{proceding-table}
\begin{array}{c|ccc}
N_+(L_{\gamma(t)}) & t<0 & t=0 & t>0\\ \hline
\Im G(i\omega)>0 & 0 & 0 & 1\\
\Im G(i\omega)<0 & 1 & 1 & 2
\end{array}
\end{align}
By complex conjugation, $N_+\!\left(L_{1/G(-i\omega)}\right)
=
N_+\!\left(L_{1/G(i\omega)}\right).$

 By Lemma~\ref{frequency}, the graphs $G(\pm i\omega)$ enclose one region if
$r\ge(p-1)/2$ and two if $r<(p-1)/2$. In the latter case,
$\Im G(i\omega)$ is negative for $0<\omega<\omega_H$ and
positive for $\omega>\omega_H$, by using  $H(0)<0$ uniqueness of
the positive zero of $H$ and \eqref{crossingderivative}.

{
For $ 1/\gamma$ strictly between the two graphs $G(\pm i\omega)$, let
$\omega>0$ satisfy $\Re \big(\frac{1}{\gamma}\big)=\Re G(i\omega)$. As shown in Table \ref{proceding-table}, one finds
\[
N_+(L_\gamma)=
\begin{cases}
0, & r\ge(p-1)/2,\\
0, & 1<r<(p-1)/2,\quad \omega>\omega_H,\\
2, & 1<r<(p-1)/2,\quad 0<\omega<\omega_H.
\end{cases}
\]
}
{If $1<r<(p-1)/2$,} at $\gamma=\gamma_H$
the only imaginary-axis eigenvalues  are $\pm i\omega_H$.
{
Taking $\frac{1}{\gamma}$ between the graphs with $\Re \big(\frac{1}{\gamma}\big)<1/\gamma_H$
and sufficiently close to $1/\gamma_H$, we have
$N_+(L_{\gamma})=0$. Lemma~\ref{localcount} then  yields
\[
0\le N_+(L_{\gamma_H})\le N_+(L_{\gamma})=0,
\]
}
which proves  part~\textup{(iii)} of the theorem.  We mention that the two conjugate arcs of \(G(\pm i\omega)\),
\(\omega>0\), that enclose the strictly stable region in the complex
\(1/\gamma\)-plane are shown in Figure~\ref{fig:stability}(d).
As exhibited in the proof, when \(1<r<(p-1)/2\), the rightmost point of
the stability boundary moves left from \(1\) to \(1/\gamma_H<1\),
so the entire strictly stable region lies in the narrower strip
\(0<\Re(1/\gamma)<1/\gamma_H\). For real \(\gamma\), strict spectral
stability holds exactly when \(\gamma>\gamma_H\) in this regime
and when \(\gamma>1\) for \((p-1)/2\le r<(p+3)/2\).

    \noindent\textbf{The case $r =\frac{p+ 3}{2}$.} 
    In this case $q=\ell$, so $f_q=f_\ell$,  the measure degenerated as
    \[
    d\mu(s)=\mu_0\,\delta_{s_0}(ds),
    \qquad
    G(\lambda)=\frac{\mu_0}{s_0-\lambda}
    \quad(\lambda\in\rho(L_0)).
    \]
    Since $G(0)=1$ as shown in Lemma \ref{moments}, we have $\mu_0=s_0$. For $\gamma\ne0$,
    in light of Lemma~\ref{rankone}, we get $s_0\in\rho(L_\gamma)$.  Therefore, substituting $G(\lambda)=s_0/(s_0-\lambda)$ into the
    characteristic equation $d_\gamma(\lambda)=1-\gamma G(\lambda)=0$, we obtain that the only possible eigenvalue in the closed right half-plane satisfying
    $\lambda=s_0(1-\gamma)$. For $\gamma\ne0$, noting that this is a simple
    zero of $d_\gamma$, so it is an algebraically simple
    eigenvalue whenever $\Re\lambda\ge0$.  For $\gamma=0$,
    the only eigenvalue of $L_0$ in the closed right half-plane
    is the simple eigenvalue $s_0$.  As a consequence, since $\Re\lambda=s_0(1-\Re\gamma)$ and $s_0>0$, we obtain
    \[
    N_+(L_\gamma)=
    \begin{cases}
    1,& \Re\gamma<1,\\
    0,& \Re\gamma\ge1.
    \end{cases}
    \]
    More precisely, when $\Re\gamma=1$, the eigenvalue
    $\lambda=s_0(1-\gamma)=-is_0\Im\gamma$ lies on the imaginary axis,
    so strict spectral stability fails.

\noindent\textbf{The case $r>(p+3)/2$.}
Since $q>\ell$, \eqref{signs} gives $\mu\ge0$, $\mu_0>0$ and $\mu((-\infty,0))>0$.$\mu$ has a finite second moment since
$L_0f_q,L_0f_{\ell+1}\in L^2(0,\infty)$.
Applying the reciprocal Cauchy transform representation (Proposition~2.2 of \cite{Maassen92})
 after normalizing and centering $\mu$, one finds
\[
\frac1{G(\lambda)}
=
\frac{\int_{\R}s\,d\mu(s)}{\mu(\R)^2}
-\frac{\lambda}{\mu(\R)}
-\int_{\R}\frac{d\nu(s)}{s-\lambda},
\]
where $\nu$ is a finite, nonzero positive measure.

By \eqref{spectral-representation-G},
$G'(x)>0$ for $0\le x<s_0$,
while $G(x)<0$ for $x>s_0$.
Since $G(0)=1$, the reciprocal $1/G$
extends holomorphically across $[0,\infty)$,
with $1/G(s_0)=0$.
By using the Stieltjes inversion formula
(Theorem~3.23 of \cite{teschl}), we obtain $\operatorname{supp}\nu\subset(-\infty,0)$.

Consequently, for $\omega>0$,
\[
\frac{d}{d\omega}\Re\frac1{G(i\omega)}
=
2\omega\int_{(-\infty,0)}
\frac{s\,d\nu(s)}{(s^2+\omega^2)^2}<0.
\]
\[
\frac1{G(i\omega)}
=-\frac{i\omega}{\mu(\mathbb R)}
+\frac{\int_{\mathbb R}s\,d\mu(s)}
{\mu(\mathbb R)^2}
+O(\omega^{-1}),
\qquad \omega\to\infty.
\]

Together with Lemma~\ref{frequency},
these properties show that
$\{1/G(i\omega):\omega\in\mathbb R\}$
is a simple proper curve intersecting
the real axis only at $\gamma=1$.
Its two conjugate branches approach the
vertical asymptote \(\Re\gamma=
\frac{\int_{\mathbb R}s\,d\mu(s)}
{\mu(\mathbb R)^2}>0.\) 
The curve therefore separates the plane
into two components.

By Lemma~\ref{rankone}, $N_+(L_\gamma)$
is constant on each component. $N_+(L_0)=1$ gives $N_+=1$
throughout the component containing $\gamma=0$.
By \eqref{spectral-representation-G} and
the preceding sign properties of $G$,
$N_+(L_\gamma)=0$ for real $\gamma>1$,
and hence in the other component.

The strict monotonicity established above
places the vertical segments between
$1/G(\pm i\omega)$, $\omega>0$,
and the region to the left of the vertical
asymptote in the component containing $0$.
The regions above and below the two boundary curves
belong to the component containing $\Re\gamma>1$.

By Lemmas~\ref{frequency} and~\ref{rankone},
together with $(1/G)'(i\omega)\ne0$,
the only imaginary-axis eigenvalue at
$\gamma=1/G(i\omega)$ is $i\omega$,
which is algebraically simple.
Applying \eqref{localcountformula} from the
component with $N_+=0$ gives
$N_+(L_\gamma)=0$ also on the curve.
This proves part~\textup{(i)}.

To prove part~\textup{(iv)}, take
$\gamma=1-\varepsilon$ with $\varepsilon>0$
sufficiently small.
For $r<(p+3)/2$, the initial count,
\eqref{simplebranch}, \eqref{criticalexpansion},
and \eqref{localcountformula} yield
\[
N_+(L_1)=
\begin{cases}
	0,&(p-1)/2\le r<(p+3)/2,\\
	1,&1<r<(p-1)/2.
\end{cases}
\]
For $r\ge(p+3)/2$, the conclusion follows
from parts~\textup{(i)} and~\textup{(ii)}. The remaining assertions in
parts~\textup{(iv)} and~\textup{(v)}
follow from Lemmas~\ref{frequency},
\ref{moments}, and~\ref{critical}.

   We remark that whenever $N_+(L_\gamma)=0$ and
    $\sigma(L_\gamma)\cap i\mathbb R=\varnothing$, we have shown that
    \[
    \sigma(L_\gamma)\cap\{\lambda\in\mathbb C:\Re\lambda\ge0\}
    =\varnothing.
    \] 
    By the spectral enclosure established in the proof
    of Lemma~\ref{rankone},
    $\sigma(L_\gamma)\cap\{\Re\lambda\ge-1/2\}$
    is contained in a compact set disjoint from
    $\sigma_{\mathrm{ess}}(L_\gamma)=(-\infty,-1]$.
    Hence, this part of the spectrum contains only
    finitely many eigenvalues. Each has a strictly
    negative real part.
    Consequently,
    \[
    \sup_{\lambda\in\sigma(L_\gamma)}
    \Re\lambda<0,
    \]
    which proves strict spectral stability by using the Definition \ref{def1-lgamma}.
    	\end{proof}

    \section*{Acknowledgment}
    The research of J.Wei is supported by National Natural Science Foundation of China (No. 12631008) and is partially supported by GRF of RGC of Hong Kong entitled “On critical and supercritical Fujita equation”. 
    
    
    \noindent {\bf AI assistance statement}: The authors acknowledge the use of AI tools. The authors reviewed and verified
    all content and results and take full responsibility
    for the  manuscript.
    \bibliographystyle{abbrv}
    \bibliography{ref}

    \end{document}